\documentclass[11pt]{amsart}
\usepackage[T1]{fontenc}
\usepackage[utf8]{inputenc}
\usepackage[english]{babel}
\usepackage{lmodern,amsmath,amssymb,amsthm,mathtools,microtype}
\usepackage{enumitem,array,booktabs}
\usepackage{tikz}
\usetikzlibrary{arrows.meta}
\usepackage[a4paper,margin=2.8cm]{geometry}
\usepackage[hypertexnames=false,hidelinks]{hyperref}
\numberwithin{equation}{section}
\makeatletter
\let\shapironormalsize\normalsize
\renewcommand{\normalsize}{\shapironormalsize
  \setlength{\abovedisplayskip}{10pt plus 2pt minus 2pt}\setlength{\belowdisplayskip}{10pt plus 2pt minus 2pt}\setlength{\abovedisplayshortskip}{7pt plus 2pt minus 1pt}\setlength{\belowdisplayshortskip}{8pt plus 2pt minus 1pt}\setlength{\jot}{4pt}}
\renewcommand{\section}{\@startsection{section}{1}{\z@}{-20pt plus 4pt minus 2pt}{10pt plus 2pt minus 1pt}{\normalfont\large\bfseries\raggedright}}
\renewcommand{\subsection}{\@startsection{subsection}{2}{\z@}{-15pt plus 3pt minus 2pt}{7pt plus 1pt minus 1pt}{\normalfont\normalsize\bfseries\raggedright}}
\renewcommand{\subsubsection}{\@startsection{subsubsection}{3}{\z@}{-12pt plus 2pt minus 1pt}{6pt plus 1pt minus 1pt}{\normalfont\normalsize\bfseries\raggedright}}
\renewcommand{\paragraph}{\@startsection{paragraph}{4}{\z@}{-10pt plus 2pt minus 1pt}{5pt plus 1pt minus 1pt}{\normalfont\normalsize\itshape\raggedright}}
\makeatother
\AtBeginDocument{\normalsize}
\theoremstyle{plain}
\newtheorem{theorem}{Theorem}[section]
\newtheorem{lemma}[theorem]{Lemma}
\newtheorem{proposition}[theorem]{Proposition}

\newtheorem{conjecture}[theorem]{Conjecture}
\newcommand{\Pp}{\mathcal P}
\newcommand{\op}[1]{\mathsf{#1}}
\newcommand{\Dd}{\mathrm{D}}
\newcommand{\duality}[2]{\langle #1,#2\rangle}
\newcommand{\ord}{\operatorname{ord}}
\newcommand{\MathTableSetup}{\renewcommand{\arraystretch}{1}\setlength{\extrarowheight}{0pt}\setlength{\aboverulesep}{5pt}\setlength{\belowrulesep}{5pt}\setlength{\abovetopsep}{2pt}\setlength{\belowbottomsep}{2pt}}
\newcommand{\MathTableStrut}{\rule[-7pt]{0pt}{25pt}}
\setlist{itemsep=7pt,topsep=7pt,parsep=3pt}
\allowdisplaybreaks[1]
\title[The Bochner--Krall problem]
{The solution of the Bochner--Krall problem}
\author{K. Castillo}
\address{CMUC, Department of Mathematics, University of Coimbra,
3000-143 Coimbra, Portugal}
\email{math@keniercastillo.com}
\author{B. Shapiro}
\address{Department of Mathematics, Stockholm University,
SE-106 91 Stockholm, Sweden, and Department of MCS,
Guangdong Technion Israel Institute of Technology,
515063 Shantou, Guangdong Province, China}
\email{shapiro@math.su.se}
\date{}
\hypersetup{pdftitle={The solution of the Bochner--Krall problem},
pdfauthor={K. Castillo and B. Shapiro}}
\begin{document}
\begin{abstract}
We solve the Bochner--Krall problem posed by Krall in 1938, extending
the classification beyond regular orthogonal sequences. We classify
monic polynomial sequences for which multiplying each polynomial by
the variable gives a linear combination of that polynomial and the
polynomials of the immediately preceding and following degrees, and
whose members are eigenpolynomials of one fixed nonscalar differential
operator of arbitrary finite order with polynomial coefficients.
Orthogonality is not assumed. The general formulation allows sequences
to start at any fixed degree. In every degree beyond a fixed threshold,
each sequence agrees, up to an invertible affine change of variable
and monic normalization, with a Hermite, Laguerre, Jacobi or Bessel
family, the monomial basis, or a finite rational Darboux transform of a
Laguerre or Jacobi family. In the positive-definite case, we determine
the representing measures, prove their uniqueness and specify the
admissible endpoint masses. No boundedness assumption on the support
or asymptotic conditions on the recurrence coefficients are imposed.
\end{abstract}
\maketitle
\enlargethispage{2pt}

\tableofcontents

\section{The Bochner--Krall problem}
\label{sec:introduction}

In 1938, H.~L. Krall asked which orthogonal polynomial sequences
consist of eigenpolynomials of one fixed nonscalar differential operator
with polynomial coefficients, of arbitrary finite order
\cite[pp.~705--706]{krall1938}. His \emph{Tchebycheff sets} are
defined by real or complex moments with nonzero Hankel determinants
in every degree. Positivity and a representing weight are not required;
the three-term recurrence follows. Krall obtained moment conditions,
proved evenness of the differential order
\cite[Theorems~2 and~3]{krall1938}, and classified the fourth-order
case in 1940 \cite{krall1940}; see also \cite[p.~600]{kwonlee2003}.

The recurrence makes multiplication by $x$ tridiagonal in the polynomial
basis; the differential equation makes a fixed differential operator
diagonal in the same basis. The problem is to classify the bases that
admit both structures, without a bound on the differential order.
It is the polynomial case of the discrete--continuous bispectral problem
\cite{gy2002,haine2001}.

The second-order equation with quadratic and linear coefficients
has antecedents in work of Euler, Gauss, Liouville and Sturm
\cite{euler1769,gauss1813,liouville1832,sturm1836}.
Routh studied its polynomial solutions, recurrences and orthogonality
in 1885 \cite{routh1885}.
Bochner's 1929 classification assumes a polynomial eigenfunction of
each degree, without assuming orthogonality \cite{bochner1929}.
His families satisfy a three-term recurrence in sufficiently large
degrees, and his list includes the canonical basis. Here the classical
orthogonal polynomial sequences are the Hermite, Laguerre, Jacobi and
Bessel families, with parameters for which their orthogonality
functionals are regular. Brenke treated integral orthogonality with
a nonnegative weight \cite{brenke1930}.

Liouville studied spectral boundary problems of arbitrary order
\cite{liouville1838}. Higher-order polynomial spectral equations
without assuming orthogonality were studied by Appell and Sheffer
\cite{appell1880montpellier,appell1880,sheffer1933}.
Orthogonality through a linear functional appears in
\cite{sheffer1931orthog}. Shohat treated generalized orthogonality
through Stieltjes integrals with an integrator of bounded variation
\cite{shohat1938}. The term \emph{Bochner--Krall}
appears in the title of A.~M. Krall and Littlejohn's 1988 report,
listed in \cite[reference~14]{ekklo1992}.

Work on the higher-order problem has followed several approaches.
Littlejohn and A.~M. Krall derived distributional equations for the
orthogonalizing functional \cite{littlejohn1984,kralllittlejohn1987}.
Koornwinder studied Jacobi weights with endpoint masses and recovered
Krall's Jacobi-type fourth-order equation by composing two second-order
operators \cite{koornwinder1984}.
For nonnegative integer density exponents, explicit finite-order
equations were constructed for Laguerre weights with a mass at zero
\cite{koekoeks1991} and symmetric Jacobi weights with equal endpoint
masses \cite{koekoek1994}. Everitt, Kwon,
Littlejohn and Wellman surveyed the classification problem and related
developments \cite{eklw2001}. For a recent account, see
\cite{shapiro2026}.

Duistermaat and Gr\"unbaum classified Schr\"odinger potentials whose
eigenfunctions also satisfy a differential equation in the spectral
parameter \cite{dg1986}.
Gr\"unbaum and Haine recovered the nonclassical families in Krall's
1940 classification by Darboux factor reversals of classical
three-term recurrences \cite{gh1996,gh1997}.
One elementary step gives the Laguerre-type family and the Jacobi-type
family with a single endpoint mass. The Legendre-type family with
equal endpoint masses is obtained in two elementary steps, although
its equation has order four \cite[p.~954]{ghh1997}.
Higher-order Laguerre and Jacobi constructions
and their algebras of differential eigenoperators followed
\cite{ghh1999,gy2002,iliev2013}; Haine surveys this development
\cite{haine2001}. The Hermite and Bessel constructions of
\cite{ghh1997} use recurrences for multiplication by a quadratic
polynomial in $x$, rather than the three-term recurrence considered here.

Wilson classified maximal rank-one bispectral algebras of ordinary
differential operators via the adelic Grassmannian \cite{wilson1993},
whose relation to Calogero--Moser spaces he developed in
\cite{wilson1998}. Here rank is the greatest common divisor of the
differential orders in the algebra. Bakalov, Horozov and Yakimov gave
an algorithmic construction of bispectral algebras of arbitrary rank
\cite{bhy1997} and related Darboux transforms of powers of Bessel
operators to highest-weight modules over $W_{1+\infty}$
\cite{bhy1998}. These differential--differential results do not impose
the three-term recurrence considered here.

In the positive case, Magnus' conjecture predicts classical densities
with possible endpoint masses \cite{kly2001}. Kwon and Lee
\cite[Theorem~5.4 and p.~618]{kwonlee2003} classified compactly
supported positive solutions as Jacobi densities with possible
endpoint masses, with integer restrictions on their exponents.
For positive measures of Nevai class whose support has convex hull
$[a,b]$, Bergkvist, Rullg\aa rd and Shapiro proved that the principal
coefficient in order $2m$ is a nonzero multiple of
$((x-a)(x-b))^m$ \cite{brs2004}.

In the matrix setting, Casper and Yakimov \cite{casperyakimov2022}
classified weight matrices admitting a symmetric second-order
differential eigenoperator under fullness of the differential-operator
algebra and positive definiteness of the leading coefficient multiplied
by the weight. Such weights arise by noncommutative bispectral Darboux
transformations from direct sums of classical weights. Their theorem
concerns matrix operators of order two, whereas the present problem is
scalar and allows arbitrary finite order.

Gr\"unbaum and Haine conjectured that all solutions of their
classification problem could be obtained from basic second-order cases
by repeated Darboux transformations \cite[pp.~145,153]{gh1996}.
Horozov, Shapiro and Tater gave an explicit list of the conjectured
families, comprising the classical families and Darboux transforms of
Laguerre and Jacobi families with the stated integer parameters
\cite[Conjecture~1.7]{hst2025}, acknowledging the earlier proposal.
That formulation retains the three-term recurrence without assuming
orthogonality. These constructions and conjectures identified the
expected families; the unresolved step was to prove exhaustiveness
in arbitrary differential order \cite[p.~2]{hst2025}.

We establish this classification in arbitrary differential order,
including Krall's regular setting and the canonical basis as a separate
alternative.
Up to an affine change and monic normalization,
every sequence satisfying the differential and three-term recurrence
conditions coincides, in all sufficiently large degrees, with a Hermite,
Laguerre, Jacobi or Bessel family, the canonical basis, or a finite rational Darboux transform of a Laguerre or Jacobi family.
This makes the Darboux construction exhaustive and extends the
classification beyond regular orthogonality.

We use Maroni's functional approach
\cite{maroni1991a,maroni1991b,maroni1993,maroni1994}:
orthogonality is represented by a linear functional on polynomials,
and dual bases translate the differential eigenvalue equations into
identities for that functional. The dual expansions are understood
in the weak dual topology, that is, by convergence on each polynomial.
See also \cite{castillopetronilho2ed}.

In the positive-definite case, the result gives a necessary and
sufficient classification of the representing measures: Hermite,
Laguerre or Jacobi densities, with possible endpoint masses in the
last two cases. A mass requires the corresponding density exponent
to be a nonnegative integer. Each measure is uniquely determined by
its moments; no boundedness assumption on the support is imposed.

\subsection{The classification conjecture}
\label{sec:statement}
Let $r\ge1$ be an integer. Consider a fixed differential operator
of order $r$ with polynomial coefficients
\begin{equation}\label{eq:differential-operator-form}
 \op{L}=\sum_{j=1}^{r}a_j(x)\partial^j,\quad \partial=\frac{d}{dx},\quad
 \deg a_j(x)\le j.
\end{equation}
The polynomial $a_r(x)$ is not identically zero, and $\deg a_j(x)=j$
for at least one $j$.
Such an operator does not increase degree. Choose $n_0\ge0$ so that, for
every $n\ge n_0$, there is a unique monic polynomial $P_n^{\op{L}}(x)$
of degree $n$, with eigenvalue $\lambda_n\in\mathbb C$, satisfying
\begin{equation}\label{eq:eigenpolynomial}
 \op{L}P_n^{\op{L}}(x)=\lambda_nP_n^{\op{L}}(x),\quad n\ge n_0.
\end{equation}
As in \cite[p.~2, following (1.2)]{hst2025}, we allow
$n_0>0$. Throughout, sufficiently large degrees means every degree
beyond one fixed integer chosen to avoid eigenvalue collisions.
Write $P(x)=(P_n(x))_{n\ge n_0}$, with $P_n(x)=P_n^{\op{L}}(x)$.
Assume that this sequence satisfies the recurrence
\begin{equation}\label{rec}
 xP_n(x)=P_{n+1}(x)+B_nP_n(x)+C_nP_{n-1}(x),\quad n>n_0,
\end{equation}
with complex coefficients $B_n,C_n$ \cite{hst2025}.

An operator $\op{M}$ is scalar if $\op{M}=c\op{I}$, $c\in\mathbb C$, where
$\op{I}$ is the identity. If $\op{M}=\sum_{j=0}^{s}b_j(x)\partial^j$ has polynomial
coefficients, with $b_s(x)$ nonzero, write $\ord \op{M}=s$.
The minimum differential order of a prescribed sequence is the least
order of a nonscalar polynomial differential operator for which every
element of that sequence is an eigenpolynomial.

\paragraph{Shifts and the two eigenvalue equations.}
For a sequence $v=(v_n)_{n\ge n_0}$, define the forward and
backward shifts by
$$
 (\op{T}v)_n=v_{n+1},\quad (\op{T}^{-1}v)_n=v_{n-1}.
$$
All compositions use only indices where the shifts are defined.
For a function $g$ of the index $n$, the operator $g(n)\op{I}$
multiplies $v_n$ by $g(n)$. We also write $g(n)$ for this operator. Products mean
composition, with the rightmost operator applied first:
$\op{A}\op{B}v=\op{A}(\op{B}v)$. In particular,
$(\op{T}(gv))_n=g(n+1)v_{n+1}$, so
$$
 \op{T}(g(n)\op{I})=g(n+1)\op{T}.
$$
Thus the recurrence operator is
$$
 \op{J}=\op{T}+B(n)\op{I}+C(n)\op{T}^{-1},\quad B(n)=B_n,\quad C(n)=C_n,
$$
and \eqref{rec} says $\op{J}P(x)=xP(x)$. Equation~\eqref{eq:eigenpolynomial}
acts on $x$, while \eqref{rec} acts on the index $n$.
A sequence satisfying both equations is called bispectral.

\paragraph{Darboux transformations.}
Let $Q(x)=(Q_n(x))_{n\ge n_0}$ be monic, with $\deg Q_n(x)=n$
and $\op{J}Q(x)=xQ(x)$. Suppose that $a\in\mathbb C$ and rational
functions $u(n),f(n)$ satisfy
$$
 B(n)-a=u(n)+f(n+1),\quad C(n)=u(n)f(n).
$$
Put
$$
 \op{U}=\op{T}+u(n)\op{I},\quad \op{F}=\op{I}+f(n)\op{T}^{-1}.
$$
Then $\op{J}-a\op{I}=\op{U}\op{F}$. Exchanging the factors defines
an elementary Darboux step:
\begin{equation}\label{eq:darboux-step}
 \widehat{\op{J}}=\op{F}\op{U}+a\op{I},\quad \widehat Q(x)=\op{F}Q(x).
\end{equation}
Each $\widehat Q_n(x)$ is monic of degree $n$. Since
$\widehat{\op{J}}\op{F}=\op{F}\op{J}$ and $\op{F}$ acts only on $n$,
$$
 \widehat{\op{J}}\widehat Q(x)=\op{F}\op{J}Q(x)=\op{F}(xQ(x))=x\widehat Q(x).
$$
These rational factor reversals are described in
\cite[Remark~1.6]{hst2025}, which also allows wider upper and lower
triangular factors. For each chain constructed below, we verify that every intermediate
sequence has a fixed nonscalar differential eigenoperator with
polynomial coefficients. The factorization of
\cite{bhy1996} is verified for the full composition in
Lemma~\ref{vf:bispectral-composition}.
In operator formulas we also write $a$ for
$a\op{I}$ and $1$ for $\op{I}$.
A finite chain has a fixed number of steps; all coefficients and
shifted polynomials are evaluated in sufficiently large degrees.
The lower bound also excludes all poles of the rational coefficients.

\paragraph{The canonical basis.}
Bochner's list also contains $P_n(x)=x^n$ \cite[p.~734]{bochner1929}.
For $n\ge0$,
\begin{equation}\label{eq:canonical-basis}
 \begin{gathered}
 xP_n(x)=P_{n+1}(x),\quad x\partial P_n(x)=nP_n(x),\\[7pt]
 (x^2\partial^2+2x\partial)P_n(x)=n(n+1)P_n(x).
 \end{gathered}
\end{equation}
The minimum nonscalar order is one. The monomials do not form a
regular orthogonal sequence.
The Bessel family is understood to have nonzero scale, and the
canonical basis is listed separately.

The following formulation is proved in
Theorem~\ref{thm:general-classification}.

\begin{conjecture}
\label{conj:final-formulation}
Let $(P_n(x))_{n\ge n_0}$ be monic, with $\deg P_n(x)=n$,
for some integer $n_0\ge0$.
Suppose that one fixed nonscalar differential operator $\op{L}$
with polynomial coefficients has every $P_n(x)$, $n\ge n_0$,
as an eigenpolynomial, and that \eqref{rec} holds for $n>n_0$.
There is an integer $N\ge n_0$ such that, after a fixed affine
change $x\mapsto ax+b$, with $a\ne0$, and monic normalization,
$(P_n(x))_{n\ge N}$ has one of the following forms:
\begin{enumerate}[label=(\roman*)]
\item The canonical basis $(x^n)_{n\ge N}$.
\item Hermite, Laguerre, Jacobi or Bessel polynomials with fixed parameters.
\item A sequence obtained by finitely many rational Darboux steps
\eqref{eq:darboux-step} from a Laguerre sequence with a positive
integer parameter or a Jacobi sequence with at least one positive
integer parameter.
\end{enumerate}
The steps and their parameters are independent of $n$. Their
identities hold for the prescribed polynomials, with their original
degree labels, for every $n\ge N$. In item~(iii), every intermediate
sequence has a fixed nonscalar differential eigenoperator with
polynomial coefficients.
No orthogonality is assumed.
\end{conjecture}

A Laguerre measure with density $e^{-x}$ on $(0,\infty)$ and a positive
mass at zero gives a one-step example; its factorization appears after
Theorem~\ref{thm:regular-hst}.

The differential eigenvalue equation makes the recurrence coefficients
eventually rational in the degree. An iterated commutator identity
and coefficient comparisons restrict the growth of the lower recurrence
coefficient to four regimes: a nonzero limit, linear growth, quadratic
growth or decay. For regular sequences, local pole and residue
restrictions in the functional identities give either a classical
family or a classical source after polynomial multiplication.
Finite connections give the rational Darboux factors; duality supplies
polynomial coefficients for every intermediate differential operator.
Without regularity, eventual rationality first separates the canonical
case. Removing and restoring a common polynomial factor handles the
remaining case while preserving the original degree labels.

Section~\ref{sec:regular-result} introduces duality and states the
regular classification, proved in Section~\ref{sec:full-regular-proof}.
Sections~\ref{sec:general-classification} and
\ref{sec:positive-classification} treat the general and positive cases.
Section~\ref{sec:concluding-remarks} records the semiclassical
consequence and the role of the recurrence assumption.

\section{The regular case}
\label{sec:regular-result}
\label{sec:functional-formulation}

\subsection{Polynomial space and duality}
The topological framework for polynomial functionals and their dual
expansions is developed in \cite[Sections~1--2]{maroni1985} and
\cite[Section~1]{maroni1988}.
Let $\Pp_k$ be the polynomials of degree at most $k$, identified
with $\mathbb C^{k+1}$ by their coefficients. Give
$$
 \Pp=\mathbb C[x]=\bigcup_{k\ge0}\Pp_k
$$
the locally convex inductive-limit topology: a linear map to a locally
convex space is continuous exactly when its restrictions to each
$\Pp_k$ are continuous. Thus all linear
functionals are continuous, as is every linear $\op{A}:\Pp\to\Pp$, since
$\op{A}(\Pp_k)$ lies in some $\Pp_\ell$.
Write $\Pp'$ for this continuous, equivalently algebraic, dual, with
$\duality{\mathbf u}{p(x)}=\mathbf u(p(x))$.

On $\Pp'$, the weak dual topology $\sigma(\Pp',\Pp)$ means
convergence on each polynomial,
or equivalently convergence of every moment
$\mu_j=\duality{\mathbf u}{x^j}$. This identifies $\Pp'$
with the product $\mathbb C^{\{0,1,\ldots\}}$, without moment
growth restrictions. The transpose $\op{A}'$, defined by
$\duality{\op{A}'\mathbf u}{p(x)}=\duality{\mathbf u}{\op{A}p(x)}$, is weakly
continuous. Multiplication and differentiation on functionals are
$$
 \begin{gathered}
 \duality{f(x)\mathbf u}{p(x)}=\duality{\mathbf u}{f(x)p(x)},\quad \Dd=-\partial',\\[7pt]
 \duality{\Dd\mathbf u}{p(x)}=-\duality{\mathbf u}{\partial p(x)},\\[7pt]
 \Dd(f(x)\mathbf u)=f(x)\Dd\mathbf u+(\partial f(x))\mathbf u.
 \end{gathered}
$$
Since $(\op{A}\op{B})'=\op{B}'\op{A}'$,
$$
 \op{L}'\mathbf u=\sum_{j=0}^{r}(-1)^j\Dd^j(a_j(x)\mathbf u)
 \quad\text{if}\quad \op{L}=\sum_{j=0}^{r}a_j(x)\partial^j.
$$
Pairings are bilinear, without complex conjugation. Point evaluation is $\boldsymbol{\delta}_a(p(x))=p(a)$;
its derivatives satisfy
$\duality{\Dd^j\boldsymbol{\delta}_a}{p(x)}=(-1)^jp^{(j)}(a)$, for $j\ge0$.

\subsection{Regularity and classification}
A functional $\mathbf u\in\Pp'$ is regular if there is a monic
sequence $(P_n(x))_{n\ge0}$, with $\deg P_n(x)=n$, such that
$$
 \duality{\mathbf u}{P_n(x)P_j(x)}=h_n\delta_{nj},\quad
 h_n\ne0\quad n,j\ge0.
$$
Here $\delta_{nj}$ is the Kronecker delta.
Since $P_0(x)=1$, we have $h_0=\mu_0\ne0$.
In particular, no nonzero polynomial $p(x)$ satisfies
$p(x)\mathbf u=0$: if $p(x)=cP_d(x)+$ lower-degree terms, then
$\duality{p(x)\mathbf u}{P_d(x)}=ch_d\ne0$.

The canonical basis cannot satisfy this definition: if $P_n(x)=x^n$,
orthogonality of $1$ and $x^2$ would give $\mathbf u(x^2)=0$,
whereas regularity would require $\mathbf u(P_1(x)^2)\ne0$.

For $n\ge0$, put $\mathbf a_n=P_n(x)\mathbf u/h_n$.
The sequence $(\mathbf a_n)_{n\ge0}$ is the dual basis to
$(P_n(x))_{n\ge0}$: for $\boldsymbol\eta\in\Pp'$,
\begin{equation}\label{eq:regular-duality}
 \duality{\mathbf a_n}{P_j(x)}=\delta_{nj},\quad
 \boldsymbol\eta=\sum_{n\ge0}\duality{\boldsymbol\eta}{P_n(x)}\mathbf a_n.
\end{equation}
The convergence in \eqref{eq:regular-duality} is in $\sigma(\Pp',\Pp)$.
For $p\in\Pp_d$, its finite expansion in $P_0,\ldots,P_d$ shows that
every partial sum through $M\ge d$ evaluates exactly to
$\duality{\boldsymbol\eta}{p(x)}$. Every fixed transpose $\op{A}'$
may be applied termwise: testing against $p(x)$ amounts to evaluating
the original expansion on the single polynomial $\op{A}p(x)$,
so again only finitely many terms contribute.
In particular, for $N\ge1$,
$\boldsymbol\eta(P_n(x))=0$ for every $n\ge N$ if and only if
$\boldsymbol\eta=q(x)\mathbf u$ for some $q(x)\in\Pp_{N-1}$
\cite[Theorems~1.3 and~6.1(i)--(ii)]{castillopetronilho2ed}.
\begin{theorem}
\label{thm:regular-hst}
Let $\mathbf u$ be regular, with monic orthogonal sequence
$P(x)=(P_n(x))_{n\ge0}$. Suppose that one fixed nonscalar differential
operator $\op{L}$, with polynomial coefficients, satisfies
$$
 \op{L}P_n(x)=\lambda_nP_n(x),\quad n\ge n_0,
$$
for some integer $n_0\ge0$.
There is an integer $N\ge n_0$ such that, after an affine change
$x\mapsto ax+b$, with $a\ne0$, and monic normalization, the sequence
$(P_n(x))_{n\ge N}$ has one of the following forms:
\begin{enumerate}[label=(\roman*)]
\item Hermite, Laguerre, Jacobi or Bessel polynomials with fixed parameters.
\item A sequence obtained by finitely many rational Darboux steps
from a monic family $Q(x)=(Q_n(x))_{n\ge n_Q}$, for some integer
$n_Q\ge0$, either Laguerre
with a positive integer parameter or Jacobi with at least one positive
integer parameter. More precisely, write $\op{J}_QQ(x)=xQ(x)$ and $\op{J}_PP(x)=xP(x)$.
There are an integer $s\ge0$, constants
$a_i\in\mathbb C$, and nonzero rational functions $f_i(n),u_i(n)$,
for $i=1,\ldots,s$, such that
\begin{equation}\label{eq:main-darboux-chain}
 \begin{gathered}
 \op{J}_{i-1}-a_i=\op{U}_i\op{F}_i,\quad \op{J}_i-a_i=\op{F}_i\op{U}_i,\\[7pt]
 \op{F}_i=1+f_i(n)\op{T}^{-1},\quad \op{U}_i=\op{T}+u_i(n),\\[7pt]
 \op{J}_0=\op{J}_Q,\quad \op{J}_s=\op{J}_P,\quad P(x)=\op{F}_s\cdots \op{F}_1Q(x).
 \end{gathered}
\end{equation}
These identities hold at every sufficiently large index, and every
intermediate sequence has a fixed nonscalar differential eigenoperator
with polynomial coefficients. For $s=0$, the product
is $\op{I}$, so $P(x)=Q(x)$.
\end{enumerate}
\end{theorem}

Regularity makes $C_n\ne0$ for $n\ge1$, so the recurrence recovers
the entire sequence backwards from any two consecutive polynomials of
sufficiently large degree.

\paragraph{A one-step example.}
The Darboux alternative of Theorem~\ref{thm:regular-hst} includes
regular orthogonal families. For $M>0$, define
the linear functionals $\mathbf u_M$ and $\mathbf v$ on
$\mathbb C[x]$ by
$$
 \mathbf u_M(p(x))=\int_0^\infty p(x)e^{-x}dx+Mp(0),\quad
 \mathbf v(p(x))=\mathbf u_M(xp(x)).
$$
The latter is the Laguerre functional with parameter one.
Write $P(x)=(P_n(x))_{n\ge0}$ and $Q(x)=(Q_n(x))_{n\ge0}$ for
their monic orthogonal sequences, and $\op{J}_P,\op{J}_Q$ for their
recurrence operators. Put $\tau_n=1+Mn$, $n\ge0$. The factors in
the Krall--Laguerre example \cite{koekoeks1991} are
$$
 \op{F}=1+\frac{n\tau_{n+1}}{\tau_n}\op{T}^{-1},\quad
 \op{U}=\op{T}+\frac{(n+1)\tau_n}{\tau_{n+1}}.
$$
For indices $n\ge1$, they satisfy
$$
 \op{J}_Q=\op{U}\op{F},\quad \op{J}_P=\op{F}\op{U},\quad
 P_n(x)=(\op{F}Q(x))_n,
$$
with $P_0(x)=Q_0(x)=1$.

\subsection{Finite products and coefficient notation}
For $k\ge0$, write
$$
 (z)_k=\prod_{j=0}^{k-1}(z+j),\quad
 z^{\underline k}=\prod_{j=0}^{k-1}(z-j).
$$
These are the rising and falling factorials; empty products equal one. We set $\deg0=-\infty$.
The notation $[t^j]p(t)$ means the coefficient of $t^j$, also for
series with negative powers. For a nonzero polynomial $p(t)$ of degree $d$,
its leading coefficient is $[t^d]p(t)$. For a differential operator
$\op{M}=\sum_j a_j(x)\partial^j$, we similarly write
$[\partial^j]\op{M}=a_j(x)$. Identities between series are read coefficient by
coefficient; convergence is specified separately.
We write $[\op{A},\op{B}]=\op{A}\op{B}-\op{B}\op{A}$ for the operator commutator.

For an operator $\op{M}$, $\mathbb C[\op{M}]$ denotes all polynomial
expressions $q(\op{M})$, with $q(z)\in\mathbb C[z]$ and $\op{M}^0=\op{I}$.

\section{Proof of the regular classification}
\label{sec:full-regular-proof}

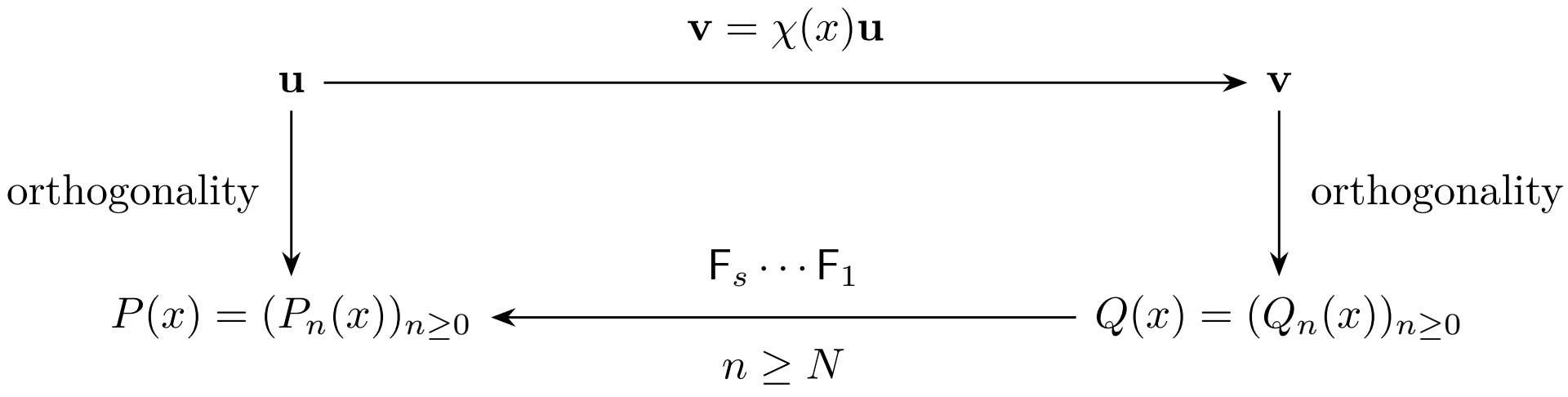
\begin{figure}[!ht]
\centering
\begin{tikzpicture}[
 font=\small,
 every node/.style={inner sep=4pt},
 >={Stealth[length=2mm,width=1.5mm]},
 line width=.55pt]
 \node (u) at (0,0) {$\mathbf u$};
 \node (v) at (8,0) {$\mathbf v$};
 \node (p) at (0,-1.9) {$P(x)=(P_n(x))_{n\ge0}$};
 \node (q) at (8,-1.9) {$Q(x)=(Q_n(x))_{n\ge0}$};
 \draw[->] (u) -- node[above=3pt]
   {$\mathbf v=\chi(x)\mathbf u$} (v);
 \draw[->] (u) -- node[left=3pt] {orthogonality} (p);
 \draw[->] (v) -- node[right=3pt] {orthogonality} (q);
 \draw[->] (q) -- node[above=3pt]
   {$\op{F}_s\cdots\op{F}_1$}
   node[below=3pt] {$n\ge N$} (p);
\end{tikzpicture}
\caption{Duality and the Darboux connection in the nonclassical regular
case, with $\mathbf v$ classical
(Lemmas~\ref{vf:classical-sources} and \ref{vf:finite-connection}).}
\label{fig:functional-darboux}
\end{figure}

We prove Lemmas~\ref{vf:foundation}--\ref{lem:bispectral-prefixes}
and then combine them to prove Theorem~\ref{thm:regular-hst}.
Figure~\ref{fig:functional-darboux} shows the nonclassical case.

\subsection{Symmetry and local restrictions}

An operator $\op{M}$ is symmetric for $\mathbf u$ when
$$
 \duality{\mathbf u}{(\op{M}p(x))q(x)}=\duality{\mathbf u}{p(x)(\op{M}q(x))}
 \quad p(x),q(x)\in\Pp.
$$

\begin{lemma}
\label{vf:foundation}
Under the hypotheses of Theorem~\ref{thm:regular-hst},
the operator $\op{L}$ does not increase degree. There is a nonconstant
polynomial $\lambda(z)\in\mathbb C[z]$ such that
$\lambda_n=\lambda(n)$ in all sufficiently large degrees.
For each fixed integer $k\ge1$, the coefficient of $x^{n-k}$ in $P_n(x)$
is rational in $n$ for sufficiently large $n$. The recurrence
$$
 xP_n(x)=P_{n+1}(x)+B_nP_n(x)+C_nP_{n-1}(x)
$$
holds in every degree, with $P_{-1}(x)=0$, $C_0=0$ and
$C_n=h_n/h_{n-1}\ne0$ for $n\ge1$. The numbers $B_n,C_n$ agree with rational functions
$B(n),C(n)$ for sufficiently large $n$.
Consecutive polynomials are coprime.

The original $\op{L}$ is symmetric for $\mathbf u$ and satisfies
$\op{L}P_n(x)=\lambda(n)P_n(x)$ for every $n\ge0$. The operator
$$
 \op{K}=\op{L}-\lambda(0)\op{I}
$$
is symmetric, has the same order $2m$, with $m\ge1$ an integer,
and satisfies $\op{K}(1)=0$. Write
$$
 \op{K}=\sum_{j=0}^{2m}k_j(x)\partial^j,\quad k_j(x)\in\Pp_j.
$$
Its leading coefficient $f(x)=k_{2m}(x)$ is a nonzero polynomial, and
\begin{equation}\label{eq:structural-equation}
 \begin{gathered}
 \Dd(f(x)\mathbf u)=\frac{k_{2m-1}(x)}m\mathbf u,\quad f(x)\Dd\mathbf u=f(x)R(x)\mathbf u,\\[7pt]
 R(x)=\frac{k_{2m-1}(x)/m-f'(x)}{f(x)},\\[7pt]
 \op{K}^*=\op{K},\quad \partial^*=-\partial-R(x).
 \end{gathered}
\end{equation}
Here $\mathbb C(x)$ denotes the rational functions of $x$.
The operation $^{*}$ is $\mathbb C$-linear, fixes every
$\varphi(x)\in\mathbb C(x)$, and satisfies
$(\op{A}\op{B})^*=\op{B}^*\op{A}^*$.
All multipliers in the functional identities, including $f(x)R(x)$,
are polynomials.
\end{lemma}

\begin{proof}
Write $\op{L}=\sum_{j=0}^r a_j(x)\partial^j$. If
$\delta=\max_j(\deg a_j(x)-j)>0$, the coefficient of
$x^{n+\delta}$ in $\op{L}P_n(x)$ is
$$
 \sum_{\deg a_j(x)-j=\delta}[x^{j+\delta}]a_j(x)
 n^{\underline j}.
$$
It vanishes for all sufficiently large integers $n$. Independence
of the falling factorials gives a contradiction. Thus
$\deg a_j(x)\le j$, and
$$
 \lambda(z)=\sum_j[x^j]a_j(x)z^{\underline j}.
$$
If $\lambda(z)$ were constant, say $c$, then $\op{L}-c\op{I}$ would
annihilate infinitely many independent polynomials and hence vanish,
contrary to nonscalarity.
Write $P_n(x)=x^n+\sum_{k=1}^n v_k(n)x^{n-k}$. For each fixed
$k\ge1$, comparison at degree $n-k$ determines $v_k(n)$
from the previous coefficients by division by
$\lambda(n)-\lambda(n-k)$. If $c$ is the leading coefficient
of $\lambda(z)$, this denominator has leading term
$ck\deg(\lambda(z))n^{\deg(\lambda(z))-1}$, so it is not the zero
polynomial. Induction proves rationality.

Orthogonality gives the stated recurrence by testing $xP_n(x)$
against $P_j(x)$. For $j<n-1$, move $x$ to the other factor;
the pairing is zero. Testing with $P_{n-1}(x)$ gives
$C_n=h_n/h_{n-1}$. In particular,
$$
 \begin{alignedat}{2}
 B_n&=v_1(n)-v_1(n+1),\quad &C_n&=v_2(n)-v_2(n+1)-B_nv_1(n)
 \end{alignedat}
$$
in large degrees, proving their rationality. The recurrence and
$C_n\ne0$ propagate a common divisor of $P_n(x),P_{n+1}(x)$
back to $P_0(x)=1$, proving coprimality.

Choose $N\ge1$ so that $\op{L}P_n(x)=\lambda(n)P_n(x)$ for $n\ge N$,
and expand by the product rule:
$$
 \begin{aligned}
 \op{A}(p(x))&=\op{L}'(p(x)\mathbf u)-(\op{L}p(x))\mathbf u\\[7pt]
 &=\sum_{j=0}^r(\partial^jp(x))\boldsymbol{\xi}_j,\quad \boldsymbol{\xi}_j\in\Pp'.
 \end{aligned}
$$
Testing against $P_j(x)$, using the eigenvalue equation for $j\ge N$
and degree preservation for $j<N$, gives
$\op{A}(P_n(x))=0$ for $n\ge N$. Since
$\duality{\op{A}(p(x))}{q(x)}=-\duality{\op{A}(q(x))}{p(x)}$,
\eqref{eq:regular-duality} gives $\op{A}(\Pp)\subset\Pp_{N-1}\mathbf u$.
Taking $p(x)=1,x,\ldots,x^r$ successively gives
$\boldsymbol{\xi}_j=f_j(x)\mathbf u$, with $f_j(x)\in\Pp$.
Thus $\sum_j f_j(x)\partial^j$ maps $\Pp$ into $\Pp_{N-1}$.
If nonzero, its largest degree shift on $x^n$ would have a nonzero
polynomial coefficient in $n$, giving outputs of unbounded degree.
It is therefore zero, so $\op{L}$ is symmetric; its triangular matrix in the orthogonal
basis $(P_n(x))_{n\ge0}$ is diagonal, giving $\op{L}P_n(x)=\lambda(n)P_n(x)$ for all $n$.
Set $\op{K}=\op{L}-\lambda(0)\op{I}$; it is nonscalar, symmetric, has the same
order as $\op{L}$, and kills $1$.

Subtracting a scalar leaves $\op{A}$ unchanged. For odd $r$, its
coefficient of $\partial^rp(x)$ would be $-2k_r(x)\mathbf u\ne0$.
Thus $r=2m>0$. The next coefficient gives
$2m\Dd(f(x)\mathbf u)-2k_{2m-1}(x)\mathbf u=0$, proving the first two
identities in \eqref{eq:structural-equation}.

Induction using the product rule and
$f(x)\Dd\mathbf u=(f(x)R(x))\mathbf u$ gives, for every polynomial $q(x)$,
$$
 \begin{gathered}
 f(x)^j\Dd^j(q(x)\mathbf u)=\bigl[f(x)^j(\partial+R(x))^jq(x)\bigr]\mathbf u,\\[7pt]
 f(x)^j(\partial+R(x))^jq(x)\in\Pp.
 \end{gathered}
$$
Substitute these identities in $\op{A}(p(x))=0$, clear powers of $f(x)$,
and use the injectivity of $q(x)\mapsto q(x)\mathbf u$.
This gives $\op{K}^*=\op{K}$ over $\mathbb C(x)$, with $\partial^*=-\partial-R(x)$.
\end{proof}

Write $\Lambda(z)=\lambda(z)-\lambda(0)$, so
$\op{K}P_n(x)=\Lambda(n)P_n(x)$ for $n\ge0$, and define the formal moment series
$$
 F(z)=\sum_{j\ge0}\duality{\mathbf u}{x^j}z^{-j-1}
 \in z^{-1}\mathbb C[[z^{-1}]].
$$
In moment series, multiplication by a polynomial $p(x)$ gives
$p(z)F(z)$ up to a polynomial, and $\Dd$ gives $F'(z)$.
Thus $\Dd(f(x)\mathbf u)=g(x)\mathbf u$, with
$g(x)=k_{2m-1}(x)/m$, gives
$$
 f(z)F'(z)+(f'(z)-g(z))F(z)\in\mathbb C[z].
$$
Dividing by $f(z)$ and using $R(z)=(g(z)-f'(z))/f(z)$, we obtain
\begin{equation}\label{eq:moment-series-equation}
 F'(z)=R(z)F(z)+S_0(z),\quad S_0(z)\in\mathbb C(z).
\end{equation}

For a nonzero rational function $\varphi(x)$,
$\ord_c\varphi(x)$ is its order of vanishing at $c$,
with negative values for poles, and $\deg_\infty\varphi(x)$ is the degree
of its numerator minus the degree
of its denominator. Its leading coefficient at infinity is the ratio
of the leading coefficients of its numerator and denominator.
We set $\deg_\infty0=-\infty$.
At a simple pole $a$ of $R(x)$, its residue is
the coefficient $\alpha$ in $R(x)=\alpha/(x-a)+O(1)$.

\begin{lemma}\label{vf:local-bounds}
Let $\op{K}$ have polynomial coefficients, order $2m>0$, and not
increase degree. Suppose $\op{K}^*=\op{K}$, where $\partial^*=-\partial-R(x)$ for
$R(x)\in\mathbb C(x)$, and write $f(x)=k_{2m}(x)$.
Then $\op{K}$ has the unique expression
$$
 \op{K}-\op{K}(1)=\sum_{j=1}^m(\partial+R(x))^jb_j(x)\partial^j,\quad
 b_j(x)\in\mathbb C(x),\quad b_m(x)=f(x).
$$
At a finite pole $R(x)\sim\rho(x-a)^{-h}$, with $h\ge2$,
$$
 \ord_a f(x)\ge mh.
$$
If $R(x)\sim\rho x^h$ at infinity, with $h\ge1$, then
$$
 \deg_\infty b_j(x)+j(h-1)\le0.
$$
If $R(x)$ tends to a nonzero constant, then
$\deg_\infty b_j(x)\le j$, and in particular $\deg f(x)\le m$.
At a simple pole with noninteger residue,
$\ord_a f(x)\ge m$.

If, in addition, $\op{K}P_n(x)=\Lambda(n)P_n(x)$ in all sufficiently large
degrees, where $(P_n(x))_{n\ge0}$ is the monic orthogonal sequence of
a regular functional $\mathbf u$, with $\deg P_n(x)=n$, and
$\Lambda(z)\in\mathbb C[z]$ is nonconstant, then $R(x)$ has no simple
pole with negative integer residue.
\end{lemma}

\begin{proof}
Degree preservation makes $\op{K}(1)$ constant.
Each summand is unchanged by $^{*}$, with leading coefficient
$b_j(x)$ at order $2j$. Subtract the highest-order term and repeat.
An odd-order remainder would change sign in its leading coefficient.
Every summand vanishes on $1$, so the remainder is $\op{K}(1)$.
Distinct highest orders give uniqueness.

Put $t=x-a$. At a pole of order $h\ge2$, write
$b_j(a+t)=b_{j,0}t^{r_j}+\cdots$, omitting zero terms.
The lowest possible exponent when the $j$-th summand acts on $t^n$
is $n+r_j-j(h+1)$, with coefficient
$b_{j,0}\rho^jn^{\underline j}$.
If the smallest degree change were $-d<-m$, its coefficient would be a
nonzero combination of falling factorials of degree at most $m$.
Since $\op{K}(t^n)$ is polynomial for $n=0,\ldots,d-1$, that
coefficient would have $d>m$ distinct roots, a contradiction.
Thus every change is at least $-m$. Taking $j=m$ gives
$\ord_a f(x)-m(h+1)\ge-m$.

At infinity with $R(x)\sim\rho x^h$, the largest degree change is
$\deg_\infty b_j(x)+j(h-1)$, with coefficient a nonzero multiple of
$n^{\underline j}$. Independence of the falling factorials rules
out a positive maximum. If $R(x)$ tends to a nonzero constant, the
same argument gives $\deg_\infty b_j(x)-j\le0$.

At a simple pole, write $R(a+t)=\alpha/t+O(1)$. Suppose that
$r_0=\ord_a f(x)<m$, and put
$$
 d=\max_{k_j\ne0}\bigl(j-\ord_a k_j(x)\bigr).
$$
The part of $\op{K}$ producing the lowest exponent is
$$
 t^{-d}p(\op{E}),\quad \op{E}=t\partial,\quad
 d\ge2m-r_0>m,\quad 0\ne p(z)\in\mathbb C[z],\quad\deg p(z)\le2m.
$$
Indeed, $c_jt^{j-d}\partial^j$ sends $t^n$ to
$c_jn^{\underline j}t^{n-d}$. The independent falling factorials
give a nonzero polynomial $p(n)$ at that power of $t$.
Since $\op{K}(t^n)$ is polynomial, $p(0)=\cdots=p(d-1)=0$.
The identities $\op{E}^*=-\op{E}-1-tR(a+t)$ and
$\op{E}t^{-d}=t^{-d}(\op{E}-d)$ show that the lowest terms in $\op{K}^*=\op{K}$ give
$$
 p(z)=p(d-\alpha-1-z).
$$
For noninteger $\alpha$, the integer roots and their reflections
are disjoint. This yields $2d>2m$ roots, a contradiction.

Finally, suppose $\alpha=-A$, with $A$ a positive integer.
Since $R(a+t)+A/t$ is analytic at zero, choose the analytic solution
$h'(t)=(R(a+t)+A/t)h(t)$, with $h(0)=1$, and define
$$
 \duality{\boldsymbol\eta}{p(x)}
 =[t^{A-1}]\bigl(h(t)p(a+t)\bigr).
$$
The coefficient of $t^{-1}$ in a derivative is zero. Since
$(t^{-A}h(t))'=R(a+t)t^{-A}h(t)$ and $\op{K}^*=\op{K}$, transferring derivatives gives
$\op{K}'\boldsymbol\eta=\op{K}(1)\boldsymbol\eta$. Hence
$$
 (\Lambda(n)-\op{K}(1))\duality{\boldsymbol\eta}{P_n(x)}=0
$$
in all sufficiently large degrees. Since $\Lambda(z)$ is nonconstant,
$\boldsymbol\eta(P_n(x))=0$ for all sufficiently large $n$.
By \eqref{eq:regular-duality}, $\boldsymbol\eta=q(x)\mathbf u$
for a polynomial $q(x)$. The identity $t^A\boldsymbol\eta=0$
then gives $t^Aq(x)\mathbf u=0$, so regularity forces $q=0$.
This contradicts $\duality{\boldsymbol\eta}{t^{A-1}}=1$.
\end{proof}

\subsection{Reduction to classical families}
\label{sec:recurrence-growth}
\label{sec:classical-sources}

Let $B(n),C(n)$ be the rational recurrence coefficients from
Lemma~\ref{vf:foundation}. Regularity gives $C_n\ne0$ for $n\ge1$.

\begin{lemma}\label{vg:recurrence-growth}
The recurrence satisfies
$$
 \deg_\infty B(n)\le 1,\quad \deg_\infty C(n)\le 2.
$$
If either $B(n)$ has degree one or $C(n)$ has degree two, then both do,
and their leading coefficients satisfy
$$
 B(n)=2\gamma n+O(1),\quad
 C(n)=\gamma^2n^2+O(n),\quad \gamma\ne0.
$$
If $C(n)$ has negative degree, then, after translating the polynomial
variable by the limit of $B(n)$,
$$
 B(n)=O(n^{-2}),\quad C(n)=O(n^{-2}).
$$
Consequently the only possible degrees of $C(n)$ are
$2,1,0$, and integers at most $-2$.
\end{lemma}

\begin{proof}
Put $r=\ord\op{K}$ and
$\op{J}=\op{T}+B(n)+C(n)\op{T}^{-1}$, with $\op{T}a(n)=a(n+1)\op{T}$.
Put $\operatorname{ad}_{\op{A}} \op{B}=[\op{A},\op{B}]$. Since
$\operatorname{ad}_x^{r+1}\op{K}=0$,
\begin{equation}\label{eq:ad-condition}
 \operatorname{ad}_{\op{J}}^{r+1}\Lambda(n)=0.
\end{equation}
The shifted polynomials have distinct degrees at large indices, so each
rational coefficient in \eqref{eq:ad-condition} vanishes at all large
integers, hence identically.

Set $b=\deg_\infty B(n)$, $p=\deg_\infty C(n)$, and
$\omega=\max\{b,p/2\}$. Suppose $\omega\ge1$. Give
$a_j(n)\op{T}^j$ weight $\deg_\infty a_j(n)+j\omega$; a finite Laurent
sum has the maximum weight of its terms. Its leading expression retains
the terms of maximum weight, replaces each coefficient by its leading
term at infinity, and replaces $\op{T}$ by $z$.
The product rule $\op{T}^ja(n)=a(n+j)\op{T}^j$ and the expansion
$a(n+j)=a(n)+ja'(n)+\cdots$ show that a commutator has weight at most
the sum of the two weights minus one. For leading expressions $\varphi(n,z),\psi(n,z)$,
its expression in that weight is
$$
 \{\varphi,\psi\}=z\bigl(\varphi_z\psi_n-\varphi_n\psi_z\bigr).
$$
Subscripts denote partial derivatives. A nonzero bracket is therefore
the leading expression; a zero bracket lowers the weight.

The leading expression of $\op{J}$ is
$$
 H(n,z)=z+\beta n^\omega+cn^{2\omega}/z.
$$
Here $\beta$ is the leading coefficient of $B(n)$ if $b=\omega$, and zero
otherwise; $c$ is the leading coefficient of $C(n)$ if $p=2\omega$, and
zero otherwise. They are not both zero. Put $\op{X}g=\{H,g\}$,
$d=\deg\Lambda(z)\ge1$, and let $\ell$ be the leading coefficient of
$\Lambda(z)$. Induction gives weight $d+k(\omega-1)$ and leading expression
$\ell\op{X}^k(n^d)$ whenever $\op{X}^k(n^d)\ne0$: the preceding iterate
is then nonzero, and the commutator rule applies, starting with
$\Lambda(n)$. Thus \eqref{eq:ad-condition} requires
$\op{X}^{r+1}(n^d)=0$. Choose a nonzero root $\xi$ of
$\xi^2+\beta\xi+c=0$, and substitute $z=\xi n^\omega$.
Then $H=0$, and differentiating the substitution agrees with
$\op{X} z$. On this branch, $\op{X}=\eta n^\omega(d/dn)$, with
$\eta=2\xi+\beta$. For half-integer $\omega$, the calculation uses powers of
$\sqrt n$. Unless $\beta^2=4c$, we have $\eta\ne0$, and

$$
 \left.\op{X}^k(n^d)\right|_{H=0}
 =\eta^k\prod_{j=0}^{k-1}\bigl(d+j(\omega-1)\bigr)
 n^{d+k(\omega-1)}\ne0.
$$
This contradicts \eqref{eq:ad-condition}.

If $\beta^2=4c$, then $p=2b$ and $\omega=b$ is an integer.
Suppose $\omega\ge2$. Set $H=1$ and put $y=\op{X} n$.
Then
$$
 y^2=1-2\beta n^\omega,\quad
 \op{X} n=y,\quad \op{X} y=-\beta\omega n^{\omega-1}.
$$
For every polynomial $F(n)$, these identities give
$$
 \op{X}^2F(n)=(1-2\beta n^\omega)F''(n)
                  -\beta\omega n^{\omega-1}F'(n).
$$
If $F(n)$ has degree $e\ge1$ and leading coefficient $\ell_F$, this
polynomial has degree $e+\omega-2$ and leading coefficient
$-\beta e(2e+\omega-2)\ell_F\ne0$. Thus every even iterate of $n^d$
is nonzero and has positive degree. A vanishing iterate would force
a later even iterate to vanish, contradicting \eqref{eq:ad-condition}.
Thus $\omega=1$, with $\beta^2=4c$.

If $p<0$, then $B(n)$ is bounded.
By Lemma~\ref{vf:foundation},
$$
 B(n)=v_1(n)-v_1(n+1),\quad
 C(n)=v_2(n)-v_2(n+1)-B(n)v_1(n),
$$
with rational $v_1(n),v_2(n)$. A rational difference has no $n^{-1}$
term at infinity. Since $B(n)$ is bounded, the polynomial part of
$v_1(n)$ has degree at most one. Translation removes its linear
part, leaving bounded $v_1(n)$ and $B(n)=O(n^{-2})$. Thus
$B(n)v_1(n)=O(n^{-2})$, and the formula for $C(n)$ has no
$n^{-1}$ term either. Its negative degree cannot be $-1$.
\end{proof}

\begin{lemma}\label{vg:classical-regimes}
With the notation of Lemma~\ref{vf:foundation}, the following four
cases exhaust the possibilities.
\begin{enumerate}[label=(\roman*),beginpenalty=10000]
\item Suppose $B(n)\to b$ and $C(n)\to c\ne0$. If
$\ord\op{K}=2m$ and $\ell$ is the leading coefficient
of $\Lambda(z)$, then
$$
 \begin{alignedat}{2}
 \deg\Lambda(z)&=2m,\quad &f(x)&=\ell\bigl((x-b)^2-4c\bigr)^m.
 \end{alignedat}
$$
The rational function $R(x)$ has the form
$$
 R(x)=\frac{\alpha}{x-a}+\frac{\beta}{x-d},\quad
 (x-a)(x-d)=(x-b)^2-4c,\quad a\ne d.
$$
Either residue may be zero.
\item Suppose $C(n)=cn+O(1)$, with $c\ne0$. Then, for a constant
$b$,
$$
 \begin{gathered}
 R(x)=-\frac{x-b}{c},\quad \op{K}\in\mathbb C[\op{H}],\\[7pt]
 \op{H}=\partial^2-\frac{x-b}{c}\partial.
 \end{gathered}
$$
The polynomials in all sufficiently large degrees are Hermite.
\item Suppose
$B(n)=2\gamma n+O(1)$ and
$C(n)=\gamma^2n^2+O(n)$, with $\gamma\ne0$.
Then
$$
 R(x)=-1/\gamma+O(1/x),\quad \deg f(x)\le m.
$$
All finite poles of $R(x)$ are simple. If one residue $\alpha$ is
noninteger, then, for its pole $a$,
$$
 \begin{gathered}
 R(x)=-1/\gamma+\frac{\alpha}{x-a},\quad \op{K}\in\mathbb C[\op{H}],\\[7pt]
 \op{H}=(x-a)\partial^2+\bigl(\alpha+1-(x-a)/\gamma\bigr)\partial.
 \end{gathered}
$$
In that case the polynomials in all large degrees are Laguerre.
\item Suppose $C(n)$ has negative degree. There are constants
$b,\alpha,\eta,\kappa$, with $\eta\ne0$ and $\kappa\ne0$, such that
$$
 \begin{gathered}
 R(x)=\frac{\alpha}{x-b}+\frac{\eta}{(x-b)^2},\\[7pt]
 f(x)=\kappa(x-b)^{2m},\quad \op{K}\in\mathbb C[\op{H}],
 \end{gathered}
$$
where
$$
 \op{H}=(x-b)^2\partial^2+\bigl((\alpha+2)(x-b)+\eta\bigr)\partial.
$$
The polynomials in all sufficiently large degrees are Bessel, and
$\deg_\infty C(n)=-2$.
\end{enumerate}
\end{lemma}

\begin{proof}
For a polynomial differential operator $\op{A}$ that does not increase
degree, put $\lambda_{\op{A}}(n)=[x^n]\op{A}(x^n)$. If $\lambda_{\op{A}}(z)$ is
nonconstant, the real polynomial $|\lambda_{\op{A}}(t)|^2$ is eventually
strictly increasing and unbounded. Thus $\lambda_{\op{A}}(n)\ne\lambda_{\op{A}}(k)$ for
$0\le k<n$ and all sufficiently large $n$. Triangularity gives
a unique monic eigenpolynomial in each such degree. Consequently
$\op{H}$ and $p(\op{H})$ have the same monic eigenpolynomials in large
degrees whenever $p(z)$ and the eigenvalue polynomial of $\op{H}$ are
nonconstant.

In cases (i)--(iii), take limits coefficient by coefficient in powers
of $z^{-1}$. Put $s_n=1,\sqrt n,n$,
respectively, and let $\op{A}_n$ be the tridiagonal $n\times n$ matrix
with entries
$$
 (\op{A}_n)_{jj}=B_{j-1}/s_n,\quad (\op{A}_n)_{j,j+1}=1,\quad
 (\op{A}_n)_{j+1,j}=C_j/s_n^2.
$$
The recurrence gives $\det(z\op{I}-\op{A}_n)=s_n^{-n}P_n(s_nz)$, hence
$$
 \begin{aligned}
 S_n(z)&:=\frac{s_nP_n'(s_nz)}{nP_n(s_nz)}\\[7pt]
 &=\sum_{k\ge0}\frac{\operatorname{tr}(\op{A}_n^k)}{nz^{k+1}}.
 \end{aligned}
$$
For fixed $k$, expand $\operatorname{tr}(\op{A}_n^k)$ over closed
paths of length $k$ \cite{bourgetmcmillen2017}. Uniformly bounded entries make the $O(k)$
boundary rows contribute $O(1/n)$ after normalization.
Closed paths have equally many upward and downward steps, so
the coefficient asymptotics give
$$
 \lim_{n\to\infty}\frac{\operatorname{tr}(\op{A}_n^k)}n
 =\begin{cases}
 [u^0](u+b+c/u)^k,&\mathrm{(i)},\\[7pt]
 0,&\mathrm{(ii)},\quad k\text{ odd},\\[7pt]
 \displaystyle\frac{c^q}{q+1}\binom{2q}{q},&\mathrm{(ii)},\quad k=2q,\\[7pt]
 \displaystyle\frac{\gamma^k}{k+1}\binom{2k}{k},&\mathrm{(iii)}.
 \end{cases}
$$
In (ii), diagonal steps vanish, leaving
$\binom{2q}{q}\int_0^1(ct)^qdt$.
In (iii), the constant term of
$(u+2\gamma t+\gamma^2t^2/u)^k$ is
$\binom{2k}{k}(\gamma t)^k$.
The coefficient asymptotics are uniform for $\varepsilon\le t\le1$.
The first $\varepsilon n$ rows contribute $O(\varepsilon)$;
letting $\varepsilon\to0$ gives the formulas for complex coefficients.
The zero distribution theorem in
\cite[Theorem~1.4]{kuijlaarsvanassche1999} additionally assumes
$B_n\in\mathbb R$ and $C_n>0$.

The binomial series gives $S_n(z)\to S(z)=z^{-1}+O(z^{-2})$, where
$$
 \begin{array}{ll}
 \mathrm{(i)}&S(z)=((z-b)^2-4c)^{-1/2},\\[7pt]
 \mathrm{(ii)}&cS(z)^2-zS(z)+1=0,\\[7pt]
 \mathrm{(iii)}&\gamma zS(z)^2-zS(z)+1=0.
 \end{array}
$$
For every fixed $j\ge0$, set
$$
 R_{j,n}(z)=\left(\frac{s_n}{n}\right)^j
 \frac{P_n^{(j)}(s_nz)}{P_n(s_nz)}.
$$
The identity
$R_{j+1,n}(z)=S_n(z)R_{j,n}(z)+n^{-1}R_{j,n}'(z)$, with $R_{0,n}(z)=1$,
gives $R_{j,n}(z)\to S(z)^j$ by induction. The series
$S_n\in z^{-1}\mathbb C[[z^{-1}]]$ and
$R_{j,n}\in z^{-j}\mathbb C[[z^{-1}]]$ have fixed lower bounds.
Thus every product coefficient is a finite sum, and multiplication
and differentiation preserve coefficientwise convergence in the
usual topology of $\mathbb C$.

\emph{(i) The bounded case.}
At $x=z$, divide the eigenvalue equation by $n^{2m}P_n(z)$.
The left side tends to $f(z)S(z)^{2m}\ne0$, so
$\deg\Lambda(z)=2m$ and $f(z)S(z)^{2m}=\ell$, giving the formula for $f(x)$. Every finite pole of $R(x)$ is a zero of $f(x)$
by its defining formula in Lemma~\ref{vf:foundation}. A pole of order
greater than one would require multiplicity at least $2m$ by
Lemma~\ref{vf:local-bounds}, whereas the two zeros of $f(x)$ each have
multiplicity $m$. Hence all poles are simple. Finally
$\deg k_{2m-1}(x)\le2m-1$ and $\deg f(x)=2m$ give $R(x)=O(1/x)$
at infinity. The partial fraction formula follows.

\emph{(ii) The Hermite case.}
Lemma~\ref{vg:recurrence-growth} makes $B(n)$ bounded. Set
$\rho=\max_{k_j\ne0}(j+\deg k_j(x))/2$, and write
$a_j=[x^{2\rho-j}]k_j(x)$ for $0\le j\le2m$.
Set $a_j=0$ for $j>2m$.
At $x=\sqrt n z$, divide the eigenvalue equation by
$n^\rho P_n(\sqrt n z)$. Its left side tends to
$$
 \sum_{j+\deg k_j(x)=2\rho} a_jz^{2\rho-j}S(z)^j.
$$
Substitute $z=(1+cS(z)^2)/S(z)$. Since $\deg k_j(x)\le j$,
we have $j\ge\rho$, and the summands become
$$
 a_j S(z)^{2j-2\rho}(1+cS(z)^2)^{2\rho-j}
$$
with distinct lowest powers. The limit is therefore a nonzero
polynomial in the nonconstant $S$, so $\rho=d=\deg\Lambda(z)$.

For these terms, $d\le j\le2d$. Put $\xi=S(z)^2$, $j=d+\nu$,
and let $\ell$ be the leading coefficient of $\Lambda(z)$.
The limiting equation is
$$
 \sum_{\nu=0}^d a_{d+\nu}\xi^\nu(1+c\xi)^{d-\nu}=\ell.
$$
Distinct lowest powers make these polynomials a basis. Expanding
$1=((1+c\xi)-c\xi)^d$ determines the coefficients. Since
$j+\deg k_j(x)\le2d$, no order above $2d$ occurs, and
$$
 \begin{gathered}
 \ord\op{K}=2d,\quad k_{2d}(x)=\ell(-c)^d,\\[7pt]
 [x]k_{2d-1}(x)=\ell d(-c)^{d-1}.
 \end{gathered}
$$
The defining formula for $R(x)$ now gives the stated linear expression.

Since $f(x)$ is constant, the structural equation gives
$\Dd\mathbf u=R(x)\mathbf u$. Hence
$\op{H}=(\partial+R(x))\partial$ is symmetric for $\mathbf u$.
Degree preservation and orthogonality give
$\op{H}P_n(x)=-(n/c)P_n(x)$ in every degree. Put
$p(z)=\Lambda(-cz)$. The operators $\op{K}$ and $p(\op{H})$ agree
on every $P_n(x)$, so $\op{K}=p(\op{H})$.
The affine change $x=b+ay$, with $a^2=2c$, identifies the normalized
polynomials with the standard Hermite eigenpolynomials.

\emph{(iii) The Laguerre case.}
Put $d=\max_j\deg k_j(x)$, let $\ell$ be the leading coefficient
of $\Lambda(z)$, and set
$A(t)=\sum_{\deg k_j(x)=d}[x^d]k_j(x)t^j$.
At $x=nz$, divide the eigenvalue equation by $n^dP_n(nz)$.
The left side tends to $z^dA(S(z))$, which is not identically zero.
The relation $z=1/(S(z)(1-\gamma S(z)))$ gives
$$
 \deg\Lambda(z)=d,\quad A(t)=\ell t^d(1-\gamma t)^d.
$$
Thus $d\ge1$ and $\ord\op{K}\ge2d$.

Suppose the polynomial part of $R(x)$ has positive degree $h$.
The infinity bound of Lemma~\ref{vf:local-bounds} gives $h=1$ and
$f(x)$ constant. The definition of $R(x)$ then makes it linear, without
finite poles. The argument in part (ii) makes $\op{K}$ a polynomial in a
Hermite operator. Its monic eigenpolynomials in sufficiently large
degrees are Hermite, contradicting $\deg_\infty C(n)=2$.
Hence $R(x)=c_0+O(1/x)$.

In $\op{K}^*$, derivatives of coefficients and terms containing $R-c_0$
lower the degree at infinity. The coefficient of $x^d$ in $\op{K}^*=\op{K}$
therefore gives
$$
 A(t)=A(-t-c_0).
$$
The reflection must exchange $0,1/\gamma$, since it has only
one fixed point. Hence $c_0=-1/\gamma\ne0$, and the infinity
bound gives $\deg f(x)\le m$. A pole of order greater than one
would require multiplicity $2m$ in $f(x)$, so none exists.

A simple pole of noninteger residue $\alpha$ requires multiplicity
at least $m$ in $f(x)$, by Lemma~\ref{vf:local-bounds}. It follows
that $f(x)=\kappa(x-a)^m$, and all finite poles of $R(x)$ must be at
$a$. This gives the asserted $R(x)$. The operator $\op{H}=(\partial+R(x))(x-a)\partial$ has polynomial
coefficients, preserves degree and is fixed by $^{*}$. At each
remaining order $2k>0$, the infinity and noninteger-residue bounds
make the leading coefficient a scalar multiple of $(x-a)^k$.
Successive subtraction gives $\op{K}\in\mathbb C[\op{H}]$.
Since $\op{H}$ has eigenvalue $-n/\gamma$, uniqueness identifies
the polynomials in large degrees as Laguerre.

\emph{(iv) The Bessel case.}
Lemma~\ref{vg:recurrence-growth} gives $B(n)\to b$ and $C(n)\to0$.
The closed-path calculation in the bounded case also applies when
$c=0$, giving $S(z)=(z-b)^{-1}$. Dividing the eigenvalue equation
by $n^{2m}P_n(z)$ as in (i) therefore gives
$$
 \deg\Lambda(z)=2m,\quad f(x)=\kappa(x-b)^{2m},\quad \kappa\ne0.
$$
Put $t=x-b$. Every finite pole of $R(x)$ is a zero of $f(x)$,
and its defining formula gives $R(x)=O(1/x)$ at infinity.
If its pole at $b$ had order $h>2$,
Lemma~\ref{vf:local-bounds} would give $mh\le2m$, a contradiction.
Consequently
$$
 R(x)=\frac{\alpha}{t}+\frac{\eta}{t^2}.
$$
If $\eta=0$, the structural equation becomes
$t^{2m}\Dd\mathbf u=\alpha t^{2m-1}\mathbf u$.
Testing it on $t^n$ gives
$$
 (n+2m+\alpha)\duality{\mathbf u}{t^{n+2m-1}}=0.
$$
The affine multiplier is nonzero for all sufficiently large integers
$n$, so the centered moments eventually vanish. Thus a power of
$t$ annihilates $\mathbf u$, contradicting regularity. Hence $\eta\ne0$.

The operator
$$
 \op{H}=t^2\partial^2+\bigl((\alpha+2)t+\eta\bigr)\partial
$$
has polynomial coefficients, preserves degree and is fixed by $*$.
At every remaining order $2j>0$, the pole bound and degree preservation
make the leading coefficient a scalar multiple of $t^{2j}$.
Successive subtraction gives $\op{K}\in\mathbb C[\op{H}]$.

The eigenvalue of $\op{H}$ in degree $n$ is $n(n+\alpha+1)$.
Uniqueness identifies $P_n(x)$ as its monic Bessel eigenpolynomials
in sufficiently large degrees.

Here $\eta$ is the Bessel scale parameter. With $b=0$ and
$\alpha=0$, its zero-scale specialization is the second-order equation
in \eqref{eq:canonical-basis} \cite[Section~11]{krallfrink1949}.
An invertible affine change cannot turn nonzero scale into zero scale,
so the canonical basis is listed separately.

In the centered expansion
$P_n(t+b)=t^n+\sum_{j=1}^n\widetilde v_j(n)t^{n-j}$,
coefficient comparison gives, for all sufficiently large $n$,
$$
\begin{gathered}
 \widetilde v_1(n)=\frac{\eta n}{2n+\alpha},\quad \widetilde v_2(n)=\frac{\eta^2n(n-1)}
 {2(2n+\alpha)(2n+\alpha-1)},\\[7pt]
 C(n)=-\frac{\eta^2n(n+\alpha)}
 {(2n+\alpha)^2(2n+\alpha-1)(2n+\alpha+1)}.
 \end{gathered}
$$
Thus $n^2C(n)\to-\eta^2/16\ne0$ and $\deg_\infty C(n)=-2$.

\end{proof}

\begin{lemma}\label{vf:classical-sources}
Consider the bounded and quadratic cases of
Lemma~\ref{vg:classical-regimes}. In the bounded case, put
$\sigma(x)=(x-a)(x-d)$. If both residues are noninteger, then
$$
 \begin{alignedat}{2}
 \op{K}&=p(\op{H}),\quad &\op{H}&=\sigma(x)\partial^2+(\sigma'(x)+\sigma(x)R(x))\partial,
 \end{alignedat}
$$
for a polynomial $p(z)$, and the $P_n(x)$ of sufficiently large degree
are Jacobi polynomials. If $C(n)$ grows quadratically and a residue
is noninteger, $\op{K}=p(\op{H})$ with $\op{H}$ as in
Lemma~\ref{vg:classical-regimes}(iii), and these polynomials are Laguerre.
In the remaining cases, the monic normalization $\chi(x)$ of $f(x)$
satisfies
\begin{equation}\label{eq:classical-functional-relation}
 \chi(x)\mathbf u=\mathbf v,
\end{equation}
where $\mathbf v$ is a regular classical Jacobi or Laguerre functional. At least one Jacobi parameter, or the Laguerre
parameter, is a positive integer.
\end{lemma}

\begin{proof}
If both residues are noninteger, the local bounds and degree
preservation make every remaining leading coefficient of order
$2j>0$ a scalar multiple of $\sigma(x)^j$. The operator $\op{H}$
in the statement has polynomial coefficients, preserves degree and is
fixed by $^{*}$. Successive subtraction gives
$\op{K}\in\mathbb C[\op{H}]$, and uniqueness identifies the
polynomials in large degrees as Jacobi.

If a residue in the bounded case is an integer, all integer residues
are nonnegative by Lemma~\ref{vf:local-bounds}.
Put
$$
 \begin{gathered}
 A=\alpha+m,\quad B=\beta+m,\\[7pt]
 \mathbf v=\sigma(x)^m\mathbf u\ne0.
 \end{gathered}
$$
The structural identity gives
$$
 \Dd(\sigma(x)\mathbf v)=((m+1)\sigma'(x)+\sigma(x)R(x))\mathbf v.
$$

After an affine change sending the endpoints to $0,1$, put
$$
 \begin{gathered}
 \phi(x)=x(1-x),\quad \psi(x)=A+1-(A+B+2)x,\\[7pt]
 \op{H}=\phi(x)\partial^2+\psi(x)\partial.
 \end{gathered}
$$
The identity $\Dd(\phi(x)\mathbf v)=\psi(x)\mathbf v$ gives
$$
 \begin{aligned}
 \duality{\mathbf v}{(\op{H}p(x))q(x)}&=-\duality{\mathbf v}{\phi(x)p'(x)q'(x)}\\[7pt]
 &=\duality{\mathbf v}{p(x)(\op{H}q(x))}.
 \end{aligned}
$$
The eigenvalues $-n(n+A+B+1)$ are pairwise distinct:
$A+B$ is either noninteger or positive. Hence $\op{H}$ has monic Jacobi
eigenpolynomials $Q_n(x)$ of every degree, pairwise orthogonal for
$\mathbf v$. A zero norm would give $Q_n(x)\mathbf v=0$ by
orthogonality, hence $\sigma(x)^mQ_n(x)\mathbf u=0$, contradicting
regularity of $\mathbf u$. Thus $\mathbf v$ is regular.

For quadratic growth, part (iii) of Lemma~\ref{vg:classical-regimes}
handles noninteger residues. Otherwise
$$
 R(x)=c+\pi'(x)/\pi(x),\quad c\ne0,
$$
where $\pi(x)$ is a nonzero polynomial and every root of $\pi(x)$ is a root of
$f(x)$. Since multiplication by $\pi(x)$ is injective on $\Pp$, its transpose is
surjective. Choose $\mathbf z$ with $\mathbf u=\pi(x)\mathbf z$.
The product rule gives
$$
 (f(x)\Dd-f(x)R(x))(\pi(x)\mathbf z)
 =f(x)\pi(x)(\Dd-c)\mathbf z=0.
$$

Put $\nabla=\partial+c$. This is invertible on each $\Pp_k$, hence on
$\Pp$; its negative transpose $\Dd-c$ is invertible on $\Pp'$.
Define
$$
 \mathbf Z_a=(\Dd-c)^{-1}\boldsymbol{\delta}_a,\quad
 \mathbf Z_a(\nabla p(x))=-p(a).
$$
A functional annihilated by a polynomial is a finite linear
combination of derivatives of point evaluations at its roots.
Dividing their polynomials in $\Dd$ by $\Dd-c$ therefore gives
$$
 \begin{gathered}
 \mathbf z=\boldsymbol\nu+\sum_a h_a\mathbf Z_a,\quad \mathbf u=\boldsymbol\rho+\sum_a h_a\pi(x)\mathbf Z_a,\\[7pt]
 \boldsymbol\rho=\pi(x)\boldsymbol\nu=\sum_a\boldsymbol\rho_a,
 \end{gathered}
$$
where $\boldsymbol\nu$ is a finite point functional supported at
roots of $f(x)$, and each $\boldsymbol\rho_a$ is supported at $a$.
Some $h_a\ne0$, since otherwise $\mathbf u$
would have a polynomial annihilator.

For each root $a$ of $f(x)$, put
$t=x-a$, $r=\ord_a f(x)\ge1$,
$\alpha=\ord_a \pi(x)$, and let $q$ be the highest derivative order in
$\boldsymbol\rho_a$, with $q=-1$ if $\boldsymbol\rho_a=0$. Since
$$
 (f(x)\Dd-f(x)R(x))(\pi(x)\mathbf Z_a)=f(a)\pi(a)\boldsymbol{\delta}_a=0,
$$
the point functional $\boldsymbol\rho$ satisfies the structural
identity. If $q\ge r-1$, its coefficient at derivative order
$q-r+1$ is a nonzero multiple of
$$
 (-1)^rf_0q^{\underline{r-1}}(q+\alpha+1),\quad
 f(a+t)=f_0t^r+O(t^{r+1}),
$$
by $t^k\Dd^j\boldsymbol{\delta}_a=(-1)^kj^{\underline k}\Dd^{j-k}\boldsymbol{\delta}_a$.
Consequently $q\le r-2$.

For polynomials $U(x),V(x)$, the identity $\op{K}^*=\op{K}$, with
$\op{K}=\sum_{k=0}^{2m}k_k(x)\partial^k$, gives
$$
 \begin{aligned}
 \nabla G(U,V)(x)&=\pi(x)\bigl((\op{K}U(x))V(x)-U(x)(\op{K}V(x))\bigr),\\[7pt]
 G(U,V)(x)&=\sum_{k=1}^{2m}\sum_{\ell=0}^{k-1}
 (-1)^\ell U^{(k-1-\ell)}(x)\nabla^\ell(\pi(x)k_k(x)V(x)).
 \end{aligned}
$$
Thus symmetry for $\mathbf u$ becomes
$$
 \boldsymbol\rho\bigl((\op{K}U(x))V(x)-U(x)(\op{K}V(x))\bigr)
 =\sum_a h_aG(U,V)(a).
$$
Interpolating the finitely many derivatives on both sides gives,
at each point,
\begin{equation}\label{eq:local-dual-symmetry}
 \boldsymbol\rho_a\bigl((\op{K}U(x))V(x)-U(x)(\op{K}V(x))\bigr)=h_aG(U,V)(a).
\end{equation}

Fix $a$ with $h_a\ne0$ and suppose $r<m$. As in
Lemma~\ref{vf:local-bounds}, the terms of maximal degree decrease are
$t^{-d}p(t\partial)$, where $N=d-\alpha-1$ and
$$
 \begin{gathered}
 d\ge2m-r>m,\quad 0\ne p(z)\in\mathbb C[z],\\[7pt]
 \deg p(z)\le2m,\\[7pt]
 p(0)=\cdots=p(d-1)=0,\quad p(z)=p(N-z).
 \end{gathered}
$$
If $\alpha\ge d$, the reflected roots give $2d>2m$ distinct
roots. Otherwise their union gives
$\alpha\le2m-d\le r$. In particular, $\alpha<d$, $N\ge1$ and $q<d$.

Write $\pi(a+t)=\pi_0t^\alpha+O(t^{\alpha+1})$ and
$p(z)=\sum_{k=d}^{2m}c_kz^{\underline k}$, with $\pi_0\ne0$.
For nonnegative integers $i,j$ with $i+j\ge N$,
$G(t^i,t^j)(a+t)=t^{i+j-N}(g(i,j)+O(t))$, where
$$
 \begin{aligned}
 g(i,j)&=\pi_0\sum_{k=d}^{2m}c_k\sum_{\ell=0}^{k-1}
 (-1)^\ell i^{\underline{k-1-\ell}}(j+k-N-1)^{\underline\ell}\\[7pt]
 &=\pi_0\frac{p(i)-p(j)}{i+j-N}\quad i+j\ne N.
 \end{aligned}
$$
Here $k_k(a+t)=c_kt^{k-d}+O(t^{k-d+1})$ for $d\le k\le2m$;
only the $\partial^\ell$ term of $\nabla^\ell$ contributes at this order.
Multiplying the sum by $i+j-N$ telescopes to
$\pi_0(p(i)-p(N-j))$; reflection gives the second line.
The quotient is polynomial, including on $i+j=N$, where
$g(i,N-i)=-\pi_0p'(N-i)=\pi_0p'(i)$.

If $q\ge0$, use $i+j=d+q$ and $0\le i\le q$ in
\eqref{eq:local-dual-symmetry}. The right side is zero, since
$i+j-N=\alpha+q+1>0$. Only the maximal degree decrease and the
highest point derivative remain on the left, so
$p(d+q-i)=p(i)=0$. Together with reflection, the roots fill
$$
 [0,d+q]\cap\mathbb Z,\quad[-\alpha-q-1,d-\alpha-1]\cap\mathbb Z,
$$
whose union has $d+\alpha+2q+2$ integers. This also holds for
$q=-1$, using the original roots and their reflections.

For integers $q+1\le i\le N-q-1$, put $j=N-i$. Then $i,j>q$,
so each product on the left of \eqref{eq:local-dual-symmetry} vanishes
to order greater than $q$, since $\op{K}$ has polynomial coefficients.
Hence $p'(i)=0$ throughout this interval, which lies in $[0,d-1]$.
Thus $\max(0,d-\alpha-2q-2)$ of the roots already counted
are double.
Consequently
$$
 \deg p(z)\ge d+\alpha+2q+2+
 \max(0,d-\alpha-2q-2)\ge2d>2m,
$$
a contradiction.

Thus every $a$ with $h_a\ne0$ has
$\ord_a f(x)\ge m$. Since $\deg f(x)\le m$, exactly
one such point remains, and
$$
 f(x)=\kappa(x-a)^m,\quad \pi(x)=\kappa_0(x-a)^\alpha,\quad
 \alpha\in\mathbb Z_{\ge0}.
$$
With $t=x-a$, put $\chi(x)=(x-a)^m$ and
$\mathbf v=\chi(x)\mathbf u\ne0$. The structural identity gives
$$
 \Dd(t\mathbf v)=(ct+\alpha+m+1)\mathbf v.
$$
The operator $\op{H}=t\partial^2+(\alpha+m+1+ct)\partial$ is symmetric for
$\mathbf v$ and has distinct eigenvalues $cn$. Its monic
eigenpolynomials $Q_n(x)$ are therefore orthogonal and, after
$y=-ct$, are Laguerre polynomials. Regularity follows as in the Jacobi case: a zero norm would give
$\chi(x)Q_n(x)\mathbf u=0$. The parameter $\alpha+m$ is a positive integer.
\end{proof}

\begin{table}[!ht]
\caption{Recurrence growth and classical sources.}
\label{tab:regimes}
\centering
\MathTableSetup
\begin{tabular}{@{}lll@{}}
\toprule
\MathTableStrut Lower coefficient & Rational function $R(x)$ & Classical source \\[7pt]
\midrule
\MathTableStrut $C(n)\to c\ne0$
 & $\dfrac{\alpha}{x-a}+\dfrac{\beta}{x-d}$ & Jacobi \\[7pt]
\MathTableStrut $C(n)\sim cn$
 & $-\dfrac{x-b}{c}$ & Hermite \\[7pt]
\MathTableStrut $C(n)\sim\gamma^2n^2$
 & $-\dfrac1\gamma+\dfrac{\alpha}{x-a}$ & Laguerre \\[7pt]
\MathTableStrut $C(n)\sim-\dfrac{\eta^2}{16n^2}$
 & $\dfrac{\alpha}{x-b}+\dfrac{\eta}{(x-b)^2}$ & Bessel \\[7pt]
\bottomrule
\end{tabular}
\end{table}

\subsection{Darboux transformations}
\label{sec:factorizations}

Finite difference operators are finite sums $\sum_i r_i(n)\op{T}^i$,
with rational coefficients acting on the left. For a nonzero such operator,
its extreme shifts are the smallest and largest indices with nonzero
coefficients; its width is the largest shift minus the smallest.

\begin{lemma}\label{lem:large-degree-duality}
Let $\mathbf u\in\Pp'$ and let $N\ge1$ be an integer.
Suppose that the monic sequence $(R_n(x))_{n\ge N}$, with
$\deg R_n(x)=n$, satisfies
$$
 \begin{aligned}
 \mathbf u(R_n(x)p(x))&=0\quad\deg p(x)<n,\\[7pt]
 \mathbf u(R_n(x)^2)&=h_n\ne0\quad n\ge N.
 \end{aligned}
$$
No nonzero polynomial annihilates $\mathbf u$.
The form $\mathbf u(p(x)q(x))$ is nondegenerate on $\Pp_{M-1}$
for every integer $M\ge N$. Every functional vanishing on $R_n(x)$,
$n\ge M$, has a unique expression $p(x)\mathbf u$, with
$\deg p(x)<M$. In particular, a functional vanishing on all sufficiently
large $R_n(x)$ and annihilated by a nonzero polynomial is zero.
Any rational differential operator having all sufficiently large
$R_n(x)$ as eigenpolynomials has polynomial coefficients.
\end{lemma}

\begin{proof}
If a monic polynomial $a(x)$ of degree $d$ annihilated $\mathbf u$,
choose $n\ge\max\{N,d\}$. Pairing $a(x)x^{n-d}-R_n(x)$ with
$R_n(x)$ gives $0=-h_n$, a contradiction.

The space $\Pp_{N-1}$ and the $R_n(x)$, $n\ge N$, form an
orthogonal decomposition of $\Pp$. If $p(x)\in\Pp_{N-1}$ is
orthogonal to $\Pp_{N-1}$, it is orthogonal to all of $\Pp$.
Thus $p(x)\mathbf u=0$, so $p=0$. Adjoining the
polynomials of nonzero norm through degree $M-1$ proves
nondegeneracy. Representing a functional on $\Pp_{M-1}$ gives the
stated representation, since both functionals vanish on the remaining $R_n(x)$.
The asserted vanishing follows from
$b(x)\boldsymbol\eta=(b(x)p(x))\mathbf u=0$ with $b(x)\ne0$, which forces $p=0$.

Let $\op{M}$ be a rational differential operator as in the statement.
At any possible finite pole $a$, every negative-power coefficient
$$
 [t^{-k}](\op{M}p)(a+t)=\boldsymbol\eta_{a,k}(p(x)),\quad k\ge1,
$$
is a finite linear combination of derivatives of point evaluation at
$a$. It vanishes on every sufficiently large $R_n(x)$ and is annihilated
by a power of $x-a$; hence $\boldsymbol\eta_{a,k}=0$ by the preceding assertion.
Thus $\op{M}$ maps every polynomial to a polynomial. Successive
evaluation on $1,x,\ldots,x^{\ord\op{M}}$ proves that each of its
coefficients is polynomial.
\end{proof}

\begin{lemma}\label{vf:finite-connection}
Let $\mathbf u,\mathbf v$ be regular, with monic orthogonal
sequences $P(x)=(P_n(x))_{n\ge0}$ and $Q(x)=(Q_n(x))_{n\ge0}$, and let
$\chi(x)\mathbf u=\mathbf v$, with $\chi(x)$ monic of degree $d$.
Suppose that, for $1\le j\le\max\{d,2\}$, the coefficients
$[x^{n-j}]P_n(x)$ and $[x^{n-j}]Q_n(x)$ agree with rational functions
of $n$ in sufficiently large degrees. Then, in sufficiently large
degrees, $P(x)$ is obtained from $Q(x)$ by finitely many rational factor reversals
$$
 \begin{gathered}
 \op{J}-a=\op{U}\op{F},\quad \op{J}_1-a=\op{F}\op{U},\\[7pt]
 Q^{(1)}(x)=\op{F}Q(x),\\[7pt]
 \op{F}=1+f(n)\op{T}^{-1},\quad \op{U}=\op{T}+u(n).
 \end{gathered}
$$
Every intermediate sequence is monic there and has a nonzero rational lower
recurrence coefficient.
\end{lemma}

\begin{proof}
Put $k_n=\mathbf v(Q_n(x)^2)\ne0$. Degree and orthogonality give
$$
 \mathbf v(P_m(x)Q_n(x))=\mathbf u(P_m(x)\chi(x)Q_n(x))=0
 \quad\text{if }m<n\text{ or }m>n+d.
$$
The dual bases give, for $n\ge d$ (see
\cite[Proposition~6.3]{castillopetronilho2ed}),
$$
 \begin{alignedat}{2}
 P_n(x)&=\sum_{j=0}^da_j(n)Q_{n-j}(x),\quad &\chi(x) Q_n(x)&=\sum_{j=0}^db_j(n)P_{n+j}(x).
 \end{alignedat}
$$
Monicity gives $a_0(n)=b_d(n)=1$, and the same pairing gives
$$
 b_j(n)h_{n+j}=a_j(n+j)k_n\quad 0\le j\le d.
$$
Thus $a_d(n)k_{n-d}=h_n$ and $b_0(n)h_n=k_n$, so both extremes
are nonzero. Comparing the first $d$ lower coefficients shows that
$a_j(n),b_j(n)$ are rational in $n$; the first two give rational
recurrence coefficients.
Regularity makes consecutive $P_n(x)$ coprime.

Let $\op{A}=1+\sum_{j=1}^da_j(n)\op{T}^{-j}$ and
$\op{B}=\sum_{j=0}^db_j(n)\op{T}^j$. Let $\op{J}_Q,\op{J}_P$ be the recurrence operators of $Q(x),P(x)$.
Then
\begin{equation}\label{eq:finite-connection}
 \begin{alignedat}{2}
 P(x)&=\op{A}Q(x),\quad &\op{A}\op{J}_Q&=\op{J}_P\op{A},\\[7pt]
 \op{B}\op{A}&=\chi(\op{J}_Q),\quad &\op{A}\op{B}&=\chi(\op{J}_P).
 \end{alignedat}
\end{equation}
These identities follow by applying both sides to $Q(x)$ or
$P(x)$, whose shifts have distinct degrees at large indices.

If $d=0$, then $\op{A}=1$. Otherwise put $\op{S}=\op{T}^{-1}$.
Scalar sequences agreeing at every sufficiently large integer index
are identified; all kernels use this convention.
In $\mathbb C(n)[\op{S}]$, with $\op{S}r(n)=r(n-1)\op{S}$,
right division leaves a remainder of smaller degree. A degree-$k$
operator with nonzero constant coefficient has a $k$-dimensional kernel,
so a remainder of degree $<k$ annihilating this kernel is zero.
The Euclidean algorithm expresses a greatest common right divisor as
a left linear combination of the two operators. Its kernel is therefore
their common kernel. In $\mathbb C(n)[\op{S}]$ it is unique up to a
nonzero rational left multiplier. After adjoining $\op{T}$, the
invertible left factors are $\rho(n)\op{T}^j$, with $\rho\ne0$ and
$j\in\mathbb Z$.
For $\deg r(x)<d$, put
$$
 \phi_r(n)=\mathbf u(r(x)Q_n(x)).
$$
Orthogonality gives $(\op{A}\phi_r)_n=\mathbf u(r(x)P_n(x))=0$
for all sufficiently large $n$.
Conversely, let $y$ satisfy $\op{A}y=0$ in sufficiently large degrees.
Choose its values at the missing indices and define $\boldsymbol\eta$
by $\boldsymbol\eta(Q_n(x))=y_n$. The identity $P(x)=\op{A}Q(x)$
gives $\boldsymbol\eta(P_n(x))=0$ in large degrees. By
\eqref{eq:regular-duality}, $\boldsymbol\eta=g(x)\mathbf u$ for a polynomial
$g(x)$. Divide $g(x)=\chi(x)g_1(x)+r(x)$, with $\deg r(x)<d$.
Then
$$
 y_n=\mathbf u(r(x)Q_n(x))+\mathbf v(g_1(x)Q_n(x))=\phi_r(n)
 \quad\text{for all sufficiently large }n.
$$
For uniqueness, if $\phi_r$ vanishes in large degrees, duality for
$Q(x)$ gives $r(x)\mathbf u=g_2(x)\mathbf v=\chi(x)g_2(x)\mathbf u$
for a polynomial $g_2(x)$. The absence of a polynomial annihilator
of $\mathbf u$ gives $r(x)=\chi(x)g_2(x)$, hence $r=0$ because
$\deg r(x)<d$. Thus every solution of $\op{A}y=0$ has the unique
form $y=\phi_r$, with $\deg r(x)<d$.

The recurrence gives $\op{J}_Q\phi_r=\phi_{\widetilde r}$, where
$\widetilde r(x)$ is the remainder of $xr(x)$ on division by $\chi(x)$.
For any root $a$ of $\chi(x)$, including a multiple root,
$\phi=\phi_{\chi(x)/(x-a)}\ne0$ therefore satisfies
$\op{A}\phi=(\op{J}_Q-a)\phi=0$.
Right division, cancelling the highest powers of $\op{S}$, gives
$$
 \op{A}=\op{R} \op{S}(\op{J}_Q-a)+r_0(n)+r_1(n)\op{S}.
$$
The remainder kills $\phi$. If it were zero,
$P(x)=(x-a)\op{R} \op{S}Q(x)$ would give a common factor of consecutive
$P_n(x)$. Neither $r_0$ nor $r_1$ is identically zero, since a nonzero single
shift cannot kill $\phi$. Normalize it to $\op{F}=1+f(n)\op{S}$, so
$\op{F}\phi=0$. Right division of $\op{S}(\op{J}_Q-a)$ by $\op{F}$
leaves a rational multiplier killing $\phi$, hence zero.
Thus $\op{F}$ divides both operators on the right, and
$$
 \op{A}=\op{A}_1\op{F},\quad \op{J}_Q-a=\op{U}\op{F},\quad \op{U}=\op{T}+u(n).
$$
Indeed, $\op{S}(\op{J}_Q-a)$ has constant coefficient $1$, so its quotient
by $\op{F}$ has the form $1+v(n)\op{S}$; multiplication by $\op{T}$ gives
the asserted $\op{U}$. Comparing the last coefficient gives
$C_Q(n)=u(n)f(n)\ne0$, where $C_Q(n)$ is the lower recurrence
coefficient of $\op{J}_Q$.

Put $\op{J}_1-a=\op{F}\op{U}$ and $Q^{(1)}(x)=\op{F}Q(x)$. Then
$$
 \op{J}_1=\op{T}+(a+u(n)+f(n))+f(n)u(n-1)\op{T}^{-1}.
$$
Its lower coefficient is nonzero, $Q^{(1)}(x)$ is monic in sufficiently
large degrees, and $P(x)=\op{A}_1Q^{(1)}(x)$, $\op{A}_1\op{J}_1=\op{J}_P\op{A}_1$.
The recursion $y_n=z_n-f(n)y_{n-1}$ solves $\op{F}y=z$.
Since $\op{A}=\op{A}_1\op{F}$, every solution of $\op{A}_1z=0$ is $z=\op{F}\phi_r$
for some $r$ of degree less than $d$. Put $\chi_1(x)=\chi(x)/(x-a)$.
The solutions of $\op{F}y=0$ are the scalar multiples of $\phi_{\chi_1(x)}$.
Dividing $r(x)$ by $\chi_1(x)$ therefore gives a unique representation
$z=\op{F}\phi_r$, with $\deg r(x)<d-1$.
Intertwining gives
$$
 \op{J}_1\op{F}\phi_r=\op{F}\phi_{\widetilde r},
$$
where $\widetilde r(x)$ is the remainder of $xr(x)$ on division by
$\chi_1(x)$. For each root $b$ of $\chi_1(x)$, the nonzero sequence
$\op{F}\phi_{\chi_1(x)/(x-b)}$ supplies the next common solution. At each step the constant coefficient is $1$,
the degree in $\op{S}$ decreases by one, and a zero remainder would
give a common factor of the final $P_n(x)$. After $d$ steps
the remaining factor is $1$, and the sequence is $P(x)$.
\end{proof}

The identity $\op{A}\op{J}_Q=\op{J}_P\op{A}$ is the Darboux
relation between recurrence operators defined in
\cite[Section~3.1, equation~(3.1)]{gy2002};
the two products in \eqref{eq:finite-connection} give the corresponding
factorizations.

\begin{lemma}\label{vf:bispectral-composition}
Let $\chi(x)\mathbf u=\mathbf v$ be the relation constructed in
Lemma~\ref{vf:classical-sources}, with $\chi(x)$ the monic normalization
of $f(x)$, and let $P(x),Q(x)$ be the corresponding monic orthogonal
sequences. After scalar normalizations of the source and transformed
sequences and, in the Jacobi case, an integer shift of the transformed
sequence, the connections \eqref{eq:finite-connection} admit the
bispectral factorization constructed in the proof below. It consists
of rational difference operators $\op{U},\op{V}$ and multiplication
by a nonzero polynomial $\theta(n)$. Each of $\op{U},\op{V},\theta$
acts on the normalized classical sequence as a rational differential
operator. The factorization and these actions are given in
\eqref{eq:bhy-factorization} and \eqref{eq:bhy-actions}, respectively.
Undoing the scalar normalizations and index shift recovers $P(x)$
with its original degree labels. The transformed differential operator
can be chosen as a nonconstant polynomial in $\op{L}$.
\end{lemma}

\begin{proof}
Normalize the Jacobi endpoints to $1,-1$ and the Laguerre endpoint
to $0$. In this proof, $\alpha,\beta$ denote the parameters of the
classical Jacobi source $\mathbf v$, and $\alpha$ denotes the parameter
of the Laguerre source. For Jacobi, choose the endpoint order so that
$\alpha$ is a positive integer, and put
$S=\alpha+\beta+1$. Put $\op{J}=\op{J}_Q$ and normalize $Q_n(x)=\kappa(n)\psi_n(x)$, where
$$
 \psi_n(x)=
 \begin{cases}
 \displaystyle\sum_{k\ge0}
 \frac{(-n)_k(n+S)_k}{(\alpha+1)_k k!}
 \left(\frac{1-x}{2}\right)^k,&\text{Jacobi},\\[7pt]
 \displaystyle\sum_{k\ge0}
 \frac{(-n)_k}{(\alpha+1)_k k!}x^k,&\text{Laguerre}.
 \end{cases}
$$
These series have rational coefficients in $n$ and terminate
for nonnegative integer $n$. The monic factors satisfy
$$
 \frac{\kappa(n+1)}{\kappa(n)}=
 \begin{cases}
 \displaystyle\frac{2(n+\alpha+1)(n+S)}{(2n+S)(2n+S+1)},
 &\text{Jacobi},\\[7pt]
 -(n+\alpha+1),&\text{Laguerre}.
 \end{cases}
$$
Replace each difference operator $\op{M}$ by $\kappa^{-1}\op{M}\kappa$,
retain its letter, and put $p_n(x)=P_n(x)/\kappa(n)$. Write
$p(x)$ and $\psi(x)$ for the corresponding columns.
The coefficients remain rational, and
$$
 \begin{alignedat}{2}
 p(x)&=\op{A}\psi(x),\quad &\op{A}\op{J}&=\op{J}_P\op{A},\\[7pt]
 \op{B}\op{A}&=\chi(\op{J}),\quad &\op{L}p_n(x)&=\lambda(n)p_n(x).
 \end{alignedat}
$$

Write $\mathcal R$ for the algebra of finite difference operators with
rational coefficients. Let $\mathcal B\subseteq\mathcal R$ consist of
those $\op{M}$ for which
\begin{equation}\label{eq:bhy-actions}
 \op{M}\psi(x)=\op{b}(\op{M})\psi(x),\quad \op{b}(\op{M})\in\mathbb C(x)[\partial].
\end{equation}
The difference operators in $n$ and rational differential operators in
$x$ act on Laurent series in $x$ for Laguerre and in $1-x$
for Jacobi, with coefficients in $\mathbb C(n)$. These actions are
defined coefficient by coefficient and commute.
A nonzero rational differential operator cannot kill polynomials of
arbitrarily large degree, so $\op{b}(\op{M})$ is unique.
Distinct shifts have distinct polynomial degrees at large integer
indices, so $\op{b}$ is injective. Commutation gives
$\op{b}(\op{M}\op{N})=\op{b}(\op{N})\op{b}(\op{M})$.
Set $\zeta(n)=n(n+S)$ for Jacobi and $\zeta(n)=n$ for
Laguerre. The differential equation and recurrence for $\psi(x)$ give
$$
 \op{b}(\op{J})=x,\quad \op{b}(\zeta)=
 \begin{cases}
 (x^2-1)\partial^2+\bigl((S+1)x+\alpha-\beta\bigr)\partial,&\text{Jacobi},\\[7pt]
 -x\partial^2+(x-\alpha-1)\partial,&\text{Laguerre}.
 \end{cases}
$$
Thus the algebra $\mathcal E$ generated by $\zeta,\op{J}$ lies in $\mathcal B$.

Identities at all large integer indices extend coefficientwise to
these series, since each coefficient is rational in $n$.

For Jacobi define the algebra automorphism $\iota$ by
$$
 \iota(r(n))=r(-n-S),\quad \iota(\op{T})=\op{T}^{-1}.
$$
Extend $\iota$ coefficientwise to $\mathbb C(n)((1-x))$, fixing $x$.
We have $\iota^2=1$, and the series gives
$\iota(\psi(x))=\psi(x)$. Applying $\iota$ to the defining identity
for $\op{M}\in\mathcal B$ shows that $\iota(\op{M})\in\mathcal B$
with the same image under $\op{b}$. Injectivity gives
$\iota(\op{M})=\op{M}$, in particular $\iota(\op{J})=\op{J}$.
Put $\op{S}_k=\op{B}\lambda(n)^k\op{A}$. Since
$$
 \begin{aligned}
 \op{S}_k\psi(x)&=\op{B} \op{L}^kp(x)\\[7pt]
 &=\op{L}^k(\op{B}p(x))=\op{L}^k(\chi(x)\psi(x)),
 \end{aligned}
$$
we have $\op{S}_k\in\mathcal B$ and $\iota(\op{S}_k)=\op{S}_k$.
Write $\op{B}=\sum_{i\in \mathcal I}b_i(n)\op{T}^i$, omitting zero coefficients.
Since $\lambda(z)$ is nonconstant, the polynomials $\lambda(n+i)$
are distinct. Vandermonde inversion expresses each $b_i(n)\op{T}^i$ as
a rational linear combination of $\op{B}\lambda^k$, $0\le k<|\mathcal I|$.
Each such monomial is invertible, giving $\op{L}_k\in\mathcal R$ with
$$
 \sum_k \op{L}_k\op{B}\lambda^k=1,\quad \op{A}=\sum_k \op{L}_k\op{S}_k.
$$
Applying $\iota$ gives $\iota(\op{A})=\op{W}\op{A}$ for some $\op{W}\in\mathcal R$.
Applying it twice and cancelling $\op{A}$ gives $\iota(\op{W})\op{W}=1$.
The extreme shifts add under products, so $\op{W}=u(n)\op{T}^j$.
Write $2s$ for the width of $\op{A}$, so $s=m$ here.
Matching supports gives $j=2s$. The argument uses only the factor
identities and $\op{L}$; for even width it requires neither equal
endpoint multiplicities nor $\deg\chi(x)=2s$. Thus, for
$\op{C}=\op{T}^s\op{A}$,
$$
 \iota(\op{C})=v\op{C},\quad v(n)=u(n-s),\quad v\iota(v)=1.
$$
Take $h=1+v$ unless $v=-1$, in which case take
$h=n+S/2$. Then $h\ne0$ and $\iota(h)v=h$. Thus
$$
 \begin{gathered}
 \begin{alignedat}{2}
 \op{A}_0&=h\op{T}^s\op{A},\quad &\op{B}_0&=\op{B}\op{T}^{-s}h^{-1},\\[7pt]
 \iota(\op{A}_0)&=\op{A}_0,\quad &\iota(\op{B}_0)&=\op{B}_0,
 \end{alignedat}\\[7pt]
 \op{B}_0\op{A}_0=\chi(\op{J}).
 \end{gathered}
$$
The assertion for $\op{B}_0$ follows by applying $\iota$ to the
last identity and cancelling $\op{A}_0$.

Every $\iota$-invariant $\op{M}\in\mathcal R$ can be multiplied on either
side by a nonzero polynomial in $\zeta$ to obtain an element of
$\mathcal E$. Suppose its extreme shifts are $-j,j$, with $j>0$,
and let $m_j(n),a_j(n)$ be the upper coefficients of
$\op{M},\op{J}^j$. Since $\mathbb C(n)$ has basis $1,n$ over
$\mathbb C(\zeta)$, write
$$
 \begin{aligned}
 \frac{m_j(n)}{a_j(n)}&=r_0(\zeta(n))+\delta_j(n)r_1(\zeta(n)),\\[7pt]
 \delta_j(n)&=\zeta(n)-\zeta(n+j)=-j(2n+S+j).
 \end{aligned}
$$
Subtract $r_0(\zeta)\op{J}^j+r_1(\zeta)[\zeta,\op{J}^j]$.
Both extreme coefficients vanish, by invariance. Induction leaves
an invariant rational function, hence a rational function of
$\zeta$. Clearing the finitely many denominators proves the left
assertion. For the right assertion, decompose the same ratio over
$\mathbb C(\zeta(n+j))$ and subtract
$\op{J}^jr_0(\zeta)+[\zeta,\op{J}^j]r_1(\zeta)$.

For Laguerre set $s=0$, $\op{A}_0=\op{A}$, and $\op{B}_0=\op{B}$.
The normalized recurrence
$\op{J}=-(n+\alpha+1)\op{T}+(2n+\alpha+1)-n\op{T}^{-1}$ gives
$$
 \begin{aligned}
 2n\op{T}^{-1}&=2n+\alpha+1-\op{J}-[n,\op{J}],\\[7pt]
 2(n+\alpha+1)\op{T}&=2n+\alpha+1-\op{J}+[n,\op{J}].
 \end{aligned}
$$
Their powers give $q_i(n)\op{T}^i\in \mathcal E$, where
$$
 q_i(n)=
 \begin{cases}
 (n+\alpha+1)_i,&i>0,\\[7pt]
 1,&i=0,\\[7pt]
 n^{\underline{-i}},&i<0.
 \end{cases}
$$
Clearing coefficient denominators and multiplying by the finitely
many shifted $q_i$'s therefore puts any $\op{M}\in\mathcal R$ in
$\mathcal E$, on either side.

In both cases choose nonzero polynomials $\mu(z),\nu(z)$ so that
\begin{equation}\label{eq:bhy-factorization}
 \begin{gathered}
 \op{U}=\mu(\zeta)\op{A}_0\in \mathcal E,\quad \op{V}=\op{B}_0\nu(\zeta)\in \mathcal E,\\[7pt]
 \theta=\nu(\zeta)\mu(\zeta),\quad \op{V}\theta^{-1}\op{U}=\chi(\op{J}),\\[7pt]
 \Psi(x)=\op{U}\psi(x)=\op{G} p(x),
 \end{gathered}
\end{equation}
where $\op{G}=\mu(\zeta)h\op{T}^s$ for Jacobi and
$\op{G}=\mu(n)$ for Laguerre. Since
$\op{U},\op{V},\theta\in\mathcal B$, $\theta$ multiplies
in $n$, and $\op{b}(\chi(\op{J}))=\chi(x)$ multiplies in $x$,
the factors satisfy the definition in \cite{bhy1996}.

Write $\theta(n)=\Theta(\zeta(n))$. For Jacobi, $\iota(\Psi(x))=\Psi(x)$
and the given equation $\op{L}\Psi(x)=\lambda(n+s)\Psi(x)$ imply
$\lambda(n+s)=\ell(\zeta(n))$ for a nonconstant polynomial
$\ell(z)$. For Laguerre take $\ell(z)=\lambda(z)$.
If $m_r$ is the multiplicity of a root $r$ of $\Theta$, put
$$
 q(w)=w\prod_{\Theta(r)=0}(w-\ell(r))^{m_r},\quad
 \rho(z)=\frac{q(\ell(z))}{\Theta(z)}\in\mathbb C[z]\setminus\{0\}.
$$
The product is over distinct roots. Divisibility holds also when
roots are multiple or values $\ell(r)$ coincide.
Replacing $\nu(z)$ by $\nu(z)\rho(z)$ multiplies $\op{V}$ on the
right by $\rho(\zeta)$, preserves the factorization, and makes
$\theta=q(\ell(\zeta))$ nonconstant.
By the bispectral Darboux theorem of \cite{bhy1996}, the factorization
\eqref{eq:bhy-factorization} gives\footnote{Alternatively, put $\op{C}_\chi=\chi(\op{J})$.
Since commutation with $\chi(x)$ lowers differential order, for
$N>\ord\op{b}(\op{U})$ the finite commutator expansion gives
$\op{C}_\chi^N\op{U}=\op{R}\op{C}_\chi$ with $\op{R}\in\mathcal B$.
The factorization, cancelling the nonzero right factor $\theta^{-1}\op{U}$,
then gives $\op{C}_\chi^N\theta=\op{R}\op{V}$. Applying $\op{b}$ to
these two identities and cancelling $\chi(x)^N$ yields
$\op{b}(\theta)=\op{b}(\op{V})\chi(x)^{-1}\op{b}(\op{U})$,
hence $\widehat{\op{L}}\Psi=\theta\Psi$.}
$$
 \begin{alignedat}{2}
 \widehat{\op{L}}&=\op{b}(\op{U})\op{b}(\op{V})\chi(x)^{-1},\quad
 &\widehat{\op{L}}\Psi(x)&=\theta(n)\Psi(x).
 \end{alignedat}
$$
Also $q(\op{L})\Psi(x)=\theta(n)\Psi(x)$. Since $\Psi_n(x)$ have arbitrarily
large degrees, $\widehat{\op{L}}=q(\op{L})$, with polynomial coefficients.
The same factors satisfy
$$
 \begin{gathered}
 \op{V}\theta^{-1}=\op{B}\op{G}^{-1},\quad \op{U}\op{V}\theta^{-1}=\op{G}\chi(\op{J}_P)\op{G}^{-1},\\[7pt]
 (\op{G} \op{J}_P\op{G}^{-1})\Psi(x)=x\Psi(x),\\[7pt]
 \op{G}^{-1}\Psi(x)=p(x),\quad P_n(x)=\kappa(n)p_n(x).
 \end{gathered}
$$
Applying $\op{G}^{-1}$ and multiplying the $n$th component by
$\kappa(n)$ recovers the prescribed $P_n(x)$ with their original
degree labels. The operator $q(\op{L})$ has polynomial coefficients
and satisfies $q(\op{L})P_n(x)=q(\lambda(n))P_n(x)$ in sufficiently
large degrees; it is nonscalar because $q\circ\lambda$ is nonconstant.
\end{proof}

\begin{lemma}\label{lem:bispectral-prefixes}
Let $\mathbf u$ and $(P_n(x))_{n\ge N}$ satisfy
Lemma~\ref{lem:large-degree-duality}, and suppose that one fixed
nonscalar differential operator with polynomial coefficients has these
polynomials as eigenpolynomials. Let
$\chi(x)\mathbf u=\mathbf v$, with $\chi(x)$ monic of degree $d$
and $\mathbf v$ a classical source from
Lemma~\ref{vf:classical-sources}. Suppose that the rational connections
\eqref{eq:finite-connection} are factored as in
Lemma~\ref{vf:finite-connection}. For any ordering of the roots of $\chi(x)$,
every intermediate sequence has a fixed nonscalar differential eigenoperator
with polynomial coefficients.
\end{lemma}

\begin{proof}
For $d=0$, the source itself gives the conclusion. Otherwise write
$$
 \begin{alignedat}{2}
 \op{A}&=\op{F}_d\cdots\op{F}_1,\quad &\op{D}_i&=\op{F}_i\cdots\op{F}_1,\\[7pt]
 P^{(i)}(x)&=\op{D}_iQ(x),\quad &\pi_i(x)&=\prod_{j=1}^i(x-a_j).
 \end{alignedat}
$$
The empty products give $\op{D}_0=1$, $P^{(0)}(x)=Q(x)$ and
$\pi_0(x)=1$. Each $\op{D}_i$ has shifts $[-i,0]$ and nonzero
extreme coefficients.

The kernel calculation in Lemma~\ref{vf:finite-connection} applies
by Lemma~\ref{lem:large-degree-duality}. With
$\phi_r(n)=\mathbf u(r(x)Q_n(x))$, each step adjoins
$\phi_{\chi/\pi_i}$ to the preceding kernel. Induction gives
$$
 \ker\op{D}_i=\{\phi_{(\chi/\pi_i)r}:\deg r<i\}.
$$
Since $\chi$ divides $\pi_i r$ exactly when $\chi/\pi_i$ divides $r$,
$$
 \ker\op{D}_i=\ker\op{A}\cap\ker\pi_i(\op{J}_Q).
$$
Right division of $\op{T}^{-i}\pi_i(\op{J}_Q)$ by $\op{D}_i$
in $\mathbb C(n)[\op{T}^{-1}]$ leaves a remainder of degree $<i$
annihilating the $i$-dimensional $\ker\op{D}_i$, hence zero.
Thus $\pi_i(\op{J}_Q)=\op{E}_i\op{D}_i$.
The greatest common right divisor of $\op{A}$ and $\pi_i(\op{J}_Q)$
has this same kernel, hence width $i$, and equals $\op{D}_i$
up to an invertible left factor.
The argument includes repeated roots.

Use the normalization in
Lemma~\ref{vf:bispectral-composition}.
For Laguerre, its denominator clearing applies to
$\op{D}_i,\op{E}_i$. It gives factors
$\op{U}=\mu(n)\op{D}_i$, $\op{V}=\op{E}_i\nu(n)$ in the source
bispectral algebra, with
$$
 \begin{alignedat}{2}
 \theta(n)&=\nu(n)\mu(n),\quad &\op{V}\theta^{-1}\op{U}&=\pi_i(\op{J}_Q).
 \end{alignedat}
$$
Multiplying $\nu$ by a further polynomial makes $\theta$
nonconstant. The same bispectral construction gives a rational differential eigenoperator for $P^{(i)}(x)$, after undoing
the scalar normalization.

For Jacobi, use its involution $\iota$. The Vandermonde argument giving
$\iota(\op{A})=w(n)\op{T}^d\op{A}$ does not require $d$ to be
even. The involution fixes $\pi_i(\op{J}_Q)$, so the description
by a common right divisor gives
$$
 \iota(\op{D}_i)=w_i(n)\op{T}^i\op{D}_i.
$$
If $i=2s$, put $v_i(n)=w_i(n-s)$ and construct $h_i$
from $v_i$ by the rule for $h$ in
Lemma~\ref{vf:bispectral-composition}, so that
$\iota(h_i)v_i=h_i$. Centering by $\op{T}^s$ and multiplying
by $h_i$ gives invariant factors whose product is $\pi_i(\op{J}_Q)$.
Two-sided denominator clearing and the same bispectral construction give a rational differential eigenoperator for
$c(n)P^{(i)}_{n+s}(x)$, with $c(n)\ne0$ for large $n$.
Choose the clearing polynomials so that $\theta$ is nonconstant.
Undoing the scalar and the integer shift gives an eigenoperator for
the prescribed $P^{(i)}_n(x)$.

If $i$ is odd, return to the monic normalization and increase the
positive integer Jacobi parameter $\alpha$ by one, at the endpoint
$1$. Denote the new monic classical sequence by $R(x)=(R_n(x))_{n\ge0}$ and put
$\mathbf v_+=(x-1)\mathbf v$. Then
$$
 \begin{gathered}
 Q(x)=\op{C}R(x),\quad \op{C}=1+c(n)\op{T}^{-1},\\[7pt]
 c(n)=-\frac{2n(n+\beta)}
 {(2n+\alpha+\beta)(2n+\alpha+\beta+1)}.
 \end{gathered}
$$
Here $c(n)$ is a nonzero rational function. Both individual parameters avoid negative integers, and their sum
remains noninteger or positive. Thus $\mathbf v_+$ is regular.
Lemma~\ref{vf:finite-connection} identifies $\op{C}$ as the single
reversal at $1$.
The connection $\op{A}\op{C}$ has multiplier $(x-1)\chi(x)$
and the corresponding kernel description. Its first $i+1$ factors have product
$\op{D}_i\op{C}$, of even width. The restriction of the recurrence
operator of $R(x)$ to $\ker(\op{D}_i\op{C})$ has characteristic polynomial
$(x-1)\pi_i(x)$. Normalize this auxiliary source as in
Lemma~\ref{vf:bispectral-composition}. Its even-width centering and
clearing give a rational differential eigenoperator for the same
$P^{(i)}_n(x)$, regardless of the endpoint multiplicities of
$(x-1)\chi(x)$.
The given sequence and its factors are unchanged.

Return to the monic normalization. To remove coefficient poles, work
backwards along the given chain.
Since $P^{(j)}(x)=\op{F}_jP^{(j-1)}(x)$, its upper factor gives
$$
 (x-a_j)P^{(j-1)}_n(x)
 =P^{(j)}_{n+1}(x)+u_j(n)P^{(j)}_n(x).
$$
Starting with $\mathbf u_d=\mathbf u$, put
$\mathbf u_{j-1}=(x-a_j)\mathbf u_j$. Degree and orthogonality give
$$
 \begin{aligned}
 \mathbf u_{j-1}(P^{(j-1)}_n(x)p(x))&=0\quad\deg p(x)<n,\\[7pt]
 \mathbf u_{j-1}(P^{(j-1)}_n(x)^2)&=u_j(n)\mathbf u_j(P^{(j)}_n(x)^2)\ne0
 \end{aligned}
$$
in large degrees. Lemma~\ref{lem:large-degree-duality} makes every rational
eigenoperator just constructed polynomial. Its eigenvalue is
nonconstant, so it is nonscalar.

\end{proof}

Under the other hypotheses of Lemma~\ref{lem:bispectral-prefixes},
the assumption that $P(x)$ already has a fixed nonscalar differential
eigenoperator with polynomial coefficients can be omitted for a
Laguerre source. The multiplier $\chi(x)$ may be any monic polynomial
for which $\chi(x)\mathbf u=\mathbf v$ and the stated rational
connections hold. Apply denominator clearing and the bispectral
construction first to $\op{A},\op{B}$, choosing $\theta$ nonconstant.
For a Jacobi source the same conclusion holds if $d=2s$ is even and,
after the conjugation by $\kappa$ in
Lemma~\ref{vf:bispectral-composition},
$$
 \iota(\op{A})=w(n)\op{T}^d\op{A},\quad w(n)\in\mathbb C(n)\setminus\{0\}.
$$
This identity gives the invariant centered factors used in that
lemma. Clearing their denominators and choosing $\theta$ nonconstant
gives a rational differential eigenoperator for $P(x)$.
In both cases, Lemma~\ref{lem:large-degree-duality} makes the
constructed eigenoperator polynomial. With this eigenoperator
constructed, the proof of Lemma~\ref{lem:bispectral-prefixes} applies
to every intermediate sequence.

\begin{proof}[Proof of Theorem~\ref{thm:regular-hst}]
Lemmas~\ref{vf:foundation}--\ref{vf:classical-sources} give either
a classical sequence in sufficiently large degrees, as in~(i), or
$\chi(x)\mathbf u=\mathbf v$, where $\mathbf v$ is a regular
Laguerre functional with a positive integer parameter or a regular
Jacobi functional with at least one positive integer parameter.
In the latter case, let $Q(x)=(Q_n(x))_{n\ge0}$ be its monic orthogonal
sequence. Lemma~\ref{vf:foundation} gives the required rationality for
both sequences, since $Q(x)$ satisfies a second-order differential equation. Lemma~\ref{vf:finite-connection} gives
\eqref{eq:main-darboux-chain}, with constants $a_i$, nonzero
rational $f_i(n),u_i(n)$, and final sequence $P(x)$.
By Lemma~\ref{vf:bispectral-composition}, the composition is
bispectral with differential operator $q(\op{L})$, $q(z)$ nonconstant,
and recovers $P(x)$ with its original indices and normalization.
Lemma~\ref{lem:bispectral-prefixes} supplies a fixed nonscalar
polynomial differential eigenoperator $\op{M}_i$ for each column
$P^{(i)}(x)$; take $\op{M}_s=\op{L}$ at the last step.
Write $\op{M}_i=\sum_jm_{ij}(x)\partial^j$.
The degree comparison in the proof of Lemma~\ref{vf:foundation},
which uses only monicity and the eigenvalue equation, gives
$\deg m_{ij}(x)\le j$ and a nonconstant eigenvalue polynomial.
Subtracting the constant zeroth-order term leaves the eigenpolynomials
unchanged and puts each operator in the form~\eqref{eq:differential-operator-form}.

Each intermediate operator has a unique monic eigenpolynomial in
every sufficiently large degree. We use the notation of
Section~\ref{sec:statement}.
The factorizations hold over $\mathbb C(n)$ at the given parameters.
There are finitely many factors and shifts, and each nonzero rational
coefficient has finitely many zeros and poles. Choose one $N$ so that,
for $n\ge N$, all identities hold, all coefficients are defined,
and every shifted degree lies in its uniqueness range.
Then, in each degree,
$$
 P_n^{(i)}(x)=P_n^{\op{M}_i}(x),\quad n\ge N.
$$
Each reversal has rational triangular factors, preserves monicity
and degree, and admits the fixed differential eigenoperator in
the form~\eqref{eq:differential-operator-form}. It is therefore a bispectral
Darboux step as defined in Section~\ref{sec:statement}.
The recurrence holds for $n>N$. The affine change is invertible,
and all parameters and the number of steps are fixed. This proves~(ii).
\end{proof}

\section{The general classification}
\label{sec:general-classification}

In the nonregular case, the named classical families are understood
through the monic eigenpolynomials of their standard second-order
operators in all sufficiently large degrees. Regular orthogonality
in all degrees is not required. The Bessel scale remains nonzero.

\begin{theorem}\label{thm:general-classification}
Let $(P_n(x))_{n\ge n_0}$ be monic, with $\deg P_n(x)=n$, and suppose
$$
 \op{L}P_n(x)=\lambda_nP_n(x)\quad n\ge n_0,
$$
for one fixed nonscalar differential operator $\op{L}$ with polynomial
coefficients. Suppose also that
$$
 xP_n(x)=P_{n+1}(x)+B_nP_n(x)+C_nP_{n-1}(x)\quad n>n_0.
$$
Then exactly one of the following alternatives holds.
\begin{enumerate}[label=(\roman*)]
\item The coefficients $C_n$ are nonzero for every sufficiently large
$n$. There is an integer $N\ge n_0$ such that
$(P_n(x))_{n\ge N}$ has one of the forms in
Theorem~\ref{thm:regular-hst}, with the same affine and monic
normalizations. In the Darboux case,
\eqref{eq:main-darboux-chain} holds for the prescribed polynomials.
\item There is $b\in\mathbb C$ such that
$$
 P_n(x)=(x-b)^n
$$
for every sufficiently large $n$. In this case $C_n=0$ there, and
the sequence restricted to sufficiently large degrees has minimum
differential order one.
\end{enumerate}
\end{theorem}

\begin{lemma}\label{lem:general-foundation}
Under the hypotheses of Theorem~\ref{thm:general-classification},
$\op{L}$ does not increase degree. There is a nonconstant polynomial
$\lambda(z)\in\mathbb C[z]$ such that $\lambda_n=\lambda(n)$
for every sufficiently large $n$. For each fixed $k\ge1$, the
coefficient $[x^{n-k}]P_n(x)$ is rational in $n$ there, as are
$B_n,C_n$.

If $C_n\ne0$ eventually, the monic greatest common divisor
$G(x)=\gcd(P_n(x),P_{n+1}(x))$ is independent of $n$
for sufficiently large $n$. Put $g=\deg G(x)$ and, for such indices,
$$
 \begin{gathered}
 q_n(x)=\frac{P_{n+g}(x)}{G(x)},\quad \lambda_0(z)=\lambda(z)-\lambda(g),\\[7pt]
 \op{L}_0=\op{L}-\lambda(g)\op{I},\quad \op{K}p(x)=\frac{\op{L}_0(G(x)p(x))}{G(x)},\\[7pt]
 \Lambda(z)=\lambda_0(z+g).
 \end{gathered}
$$
The operator $\op{K}$ has polynomial coefficients, does not increase degree,
annihilates $1$, and satisfies $\op{K}q_n(x)=\Lambda(n)q_n(x)$ in sufficiently large
degrees. Consecutive $q_n(x)$ are coprime. There are an integer
$N\ge1$ and a functional $\mathbf u$ such that
\begin{equation}\label{eq:general-orthogonality}
 \mathbf u(q_n(x)p(x))=0\quad\deg p(x)<n,\quad
 \mathbf u(q_n(x)^2)=h_n\ne0\quad n\ge N.
\end{equation}
No nonzero polynomial annihilates $\mathbf u$, and the form
$\mathbf u(p(x)q(x))$ is nondegenerate on $\Pp_{M-1}$ for every
$M\ge N$. Every functional vanishing on all $q_n(x)$, $n\ge M$,
therefore has a unique expression $r(x)\mathbf u$, with $\deg r(x)<M$.
The operator $\op{K}$ is symmetric for $\mathbf u$ and has positive
even order $2m$, with $m\ge1$. Write
$$
 \begin{gathered}
 \op{K}=\sum_{j=0}^{2m}k_j(x)\partial^j,\quad f(x)=k_{2m}(x),\\[7pt]
 R(x)=\frac{k_{2m-1}(x)/m-f'(x)}{f(x)}.
 \end{gathered}
$$
With this notation, \eqref{eq:structural-equation} holds.
\end{lemma}

\begin{proof}
Write $\op{L}=\sum_j a_j(x)\partial^j$.
The degree and coefficient comparisons in the proof of
Lemma~\ref{vf:foundation} use only monicity and the differential
equation. They give $\deg a_j(x)\le j$, a nonconstant eigenvalue
polynomial $\lambda(z)$, and rational coefficients
$[x^{n-k}]P_n(x)$ for each fixed $k$. The formulas for $B_n,C_n$
in that proof give their rationality as well.

When $C_n\ne0$, the recurrence gives the equality of ideals
$(P_n(x),P_{n+1}(x))=(P_{n-1}(x),P_n(x))$. Thus $G(x)$ is independent of $n$ in large degrees. Dividing by it
gives coprime polynomials with recurrence coefficients $B(n+g),C(n+g)$.
The operator $\op{K}$ initially has rational coefficients and satisfies
$\op{K}q_n(x)=\Lambda(n)q_n(x)$ in sufficiently large degrees.

Rename the recurrence coefficients of $q_n(x)$ as $B_n,C_n$, and
choose $N\ge1$ so that all preceding identities hold for $n\ge N$. B\'ezout's identity gives
polynomials $s_N(x),s_{N+1}(x)$, with degrees below $N,N+1$, such that
$$
 q_N(x)s_{N+1}(x)-q_{N+1}(x)s_N(x)=h_N\ne0.
$$
Propagate $s_n(x)$ by the same recurrence as $q_n(x)$. Then
$\deg s_n(x)<n$, and direct subtraction gives
$$
 \begin{gathered}
 q_n(x)s_{n+1}(x)-q_{n+1}(x)s_n(x)=h_n,\quad h_n=C_nh_{n-1}\ne0\quad n>N,\\[7pt]
 \frac{s_{n+1}(z)}{q_{n+1}(z)}-\frac{s_n(z)}{q_n(z)}=\frac{h_n}{q_n(z)q_{n+1}(z)}\quad n\ge N.
 \end{gathered}
$$
The last expression is $O(z^{-2n-1})$ at infinity. Thus the
coefficients of $s_n(z)/q_n(z)$ stabilize and define a series
$F(z)=\sum_{j\ge0}\mu_jz^{-j-1}$. Summing the displayed differences
gives
\begin{equation}\label{eq:general-moment-approximation}
 q_n(z)F(z)-s_n(z)=h_nz^{-n-1}+O(z^{-n-2}).
\end{equation}
Set $\mathbf u(x^j)=\mu_j$. The coefficient of $z^{-k-1}$ in
$q_n(z)F(z)$ is $\mathbf u(q_n(x)x^k)$, proving
\eqref{eq:general-orthogonality}.

Lemma~\ref{lem:large-degree-duality} gives the absence of a polynomial
annihilator, the stated nondegeneracy and dual representation, and
polynomial coefficients for $\op{K}$. Writing $\op{K}=\sum_jk_j(x)\partial^j$, degree
comparison gives $\deg k_j(x)\le j$. Its diagonal polynomial is
$\Lambda(z)$, since they agree at all sufficiently large
integers, so $\op{K}(1)=\Lambda(0)=0$.

Put $\mathcal W=\Pp_{N-1}$. It is $\op{K}$-invariant and orthogonal
to every $q_n(x)$, $n\ge N$. Hence
$\op{K}'(q_n(x)\mathbf u)=\Lambda(n)q_n(x)\mathbf u$ for $n\ge N$.
The symmetry argument in the proof of Lemma~\ref{vf:foundation}
therefore applies, using Lemma~\ref{lem:large-degree-duality} in place
of \eqref{eq:regular-duality}. Its coefficient comparisons give
symmetry, positive even order, and \eqref{eq:structural-equation}.
\end{proof}

\begin{lemma}\label{lem:common-factor-bound}
Use $G,\op{K},\mathbf u$ from Lemma~\ref{lem:general-foundation}.
If, at a root $a$ of $G(x)$,
$$
 R(a+t)=\frac{\alpha}{t}+R_0(t),\quad
 \alpha\in\mathbb Z_{\ge0},
$$
where $R_0(t)$ is analytic at zero, then
$\ord_aG(x)\le\alpha$.
\end{lemma}

\begin{proof}
Put $g_a=\ord_aG(x)$. Polynomiality of
$G(x)\op{K}G(x)^{-1}$ shows that $\op{K}$ preserves local functions
with a pole of order at most $g_a$.
Choose $h(t)$ analytic and nonzero near zero, with $h'(t)/h(t)=R_0(t)$,
and set $w(t)=t^\alpha h(t)$. For those local functions $f(t)$,
define
$$
 \boldsymbol\eta_f(p(x))
 =\operatorname{Res}_{t=0}w(t)f(t)p(a+t).
$$
Here the residue is the coefficient of $t^{-1}$.
Their image is the space of point functionals of derivative orders
at most $g_a-\alpha-1$, of dimension $\max(g_a-\alpha,0)$.
Since $\op{K}^*=\op{K}$, integration by parts in the residue gives
$\op{K}'\boldsymbol\eta_f=\boldsymbol\eta_{\op{K}f}$.
If $g_a>\alpha$, this nonzero finite-dimensional invariant space
contains an eigenfunctional $\boldsymbol\eta$, with
$\op{K}'\boldsymbol\eta=c\boldsymbol\eta$. Thus
$$
 (\Lambda(n)-c)\boldsymbol\eta(q_n(x))=0.
$$
Since $\Lambda(z)$ is nonconstant and $\boldsymbol\eta$ has a point
annihilator, Lemma~\ref{lem:large-degree-duality} gives $\boldsymbol\eta=0$,
a contradiction.
\end{proof}

\begin{lemma}\label{lem:polynomial-conjugates}
Let $\op{M}$ be a polynomial differential operator that does not
increase degree and satisfies $\op{M}^*=\op{M}$, with
$\partial^*=-\partial-R_M(x)$.
Let $H(x)$ be a nonconstant polynomial whose distinct roots $a_i$
have multiplicities at most $\alpha_i$, where
$$
 R_M(x)=\frac{\alpha_i}{x-a_i}+O(1),\quad
 \alpha_i\in\mathbb Z_{\ge1}.
$$
There is a nonconstant polynomial $c(z)$ such that
$H(x)c(\op{M})H(x)^{-1}$ has polynomial coefficients and does not
increase degree. If, for every sufficiently large $n$, there is a
monic polynomial $p_n(x)$ of degree $n$ such that
$\op{M}p_n(x)=\lambda_M(n)p_n(x)$, where
$\lambda_M(z)\in\mathbb C[z]$ is nonconstant, this conjugate is nonscalar.
\end{lemma}

\begin{proof}
At $a_i$, put $t=x-a_i$, and choose an analytic weight
$w_i(t)=t^{\alpha_i}h_i(t)$, with $h_i(t)$ nonvanishing and
$w_i'(t)/w_i(t)=R_M(a_i+t)$. If
$\op{M}=\sum_jm_j(x)\partial^j$, weighted symmetry gives
$$
 w_i\op{M}f=\sum_j(-1)^j\partial^j(m_jw_if).
$$
If $f(x)$ has a pole of order at most $\alpha_i$, the right side
is analytic. Hence $\op{M}f$ has the same pole bound.
Let $\mathcal V$ be the rational functions whose only poles are
at the $a_i$, with these bounds. The operator preserves both
$\mathcal V$ and $\Pp$, and the quotient $\mathcal V/\Pp$
has dimension $\sum_i\alpha_i$. Let $c(z)$ be the characteristic
polynomial of its induced map. Then $c(\op{M})\mathcal V\subset\Pp$.
In particular, $H(x)c(\op{M})H(x)^{-1}$ maps every polynomial to a
polynomial. Evaluation on successive monomials shows that its
coefficients are polynomial. Since $\deg m_j(x)\le j$, the operator
does not increase degree at infinity even on rational functions;
conjugation by $H(x)$ preserves this bound.
If $\op{M}p_n(x)=\lambda_M(n)p_n(x)$ with a nonconstant polynomial
$\lambda_M(z)$, the conjugate has nonconstant eigenvalue
$c(\lambda_M(n))$ on $H(x)p_n(x)$.
\end{proof}

\begin{proof}[Proof of Theorem~\ref{thm:general-classification}]
By Lemma~\ref{lem:general-foundation}, either $C(n)=0$ or
$C_n\ne0$ for every sufficiently large $n$.

\emph{The vanishing case.}
Choose $N$ so that $P_{n+1}(x)=(x-B_n)P_n(x)$ and the differential
equation hold for $n\ge N$. Its order $r$ is positive: a
polynomial multiplier having a nonzero eigenpolynomial is scalar.
For every $k\ge N$, the quotients
$P_n(x)/P_k(x)$, $n\ge k$, form a monic basis of $\Pp$.
Consequently $\op{L}$ preserves $P_k(x)\Pp$, and
$P_k(x)^{-1}\op{L}P_k(x)$ has polynomial coefficients, as follows
by testing $1,x,\ldots,x^r$. Its coefficient of $\partial^{r-1}$ is
$$
 a_{r-1}(x)+ra_r(x)\frac{P_k'(x)}{P_k(x)}.
$$
For $k\ge N$, every zero of $P_k(x)$ is therefore a zero of $a_r(x)$.
Since $B_n$ is a zero of $P_{n+1}(x)$, the rational function $B(n)$
takes values in the fixed finite zero set of $a_r(x)$ at all large
integers. Hence $B(n)=b$ for a constant $b$, and
$$
 P_n(x)=P_N(x)(x-b)^{n-N}\quad n\ge N.
$$
The polynomial operator $P_N(x)^{-1}\op{L}P_N(x)$ diagonalizes every
power $(x-b)^k$, $k\ge0$. Comparing coefficients on these powers
successively shows that its coefficient of $\partial^j$ is
$c_j(x-b)^j$, with $c_j\in\mathbb C$. Its leading coefficient is
unchanged by conjugation, so $a_r(x)=c_r(x-b)^r$, with $c_r\ne0$.
The established root condition gives $P_N(x)=(x-b)^N$, hence
$P_n(x)=(x-b)^n$ for $n\ge N$. The operator $(x-b)\partial$ has order one and diagonalizes
these polynomials, proving~(ii).

\emph{The nonvanishing case.}
Use $G,q_n,\op{K},\mathbf u$ from
Lemma~\ref{lem:general-foundation}. Write $\ord\op{K}=2m$, choose
$N$ as in that lemma, and put $\mathcal W=\Pp_{N-1}$. In the remainder
of the proof, write $B_n,C_n$ for the recurrence coefficients of
$q_n(x)$. Lemmas~\ref{vf:local-bounds} and \ref{vg:recurrence-growth} apply
with these degree, symmetry and rationality properties. In the exclusion
of a negative integer residue, only the tail duality of
Lemma~\ref{lem:large-degree-duality} is needed.

For all four regimes of Lemma~\ref{vg:classical-regimes},
complete the missing degrees by arbitrary monic polynomials and let
$\op{M}_n$ be the multiplication matrix modulo $q_n(x)$. Its
characteristic polynomial is $q_n(z)$; outside a fixed initial block
it is the tridiagonal recurrence matrix. Take $s=1$ in the bounded
and decreasing cases, $s=\sqrt n$ in the linear-growth case, and
$s=n$ in the quadratic case. Diagonal conjugation of $\op{M}_n/s$
multiplies an entry $(i,j)$ of the initial block by $s^{j-i-1}$.
Since $j\le i+1$, these factors do not grow. The row sums remain
uniformly bounded, and a closed path of fixed length visiting the
initial block contributes only $O(1)$ to the trace. The normalized
trace limits and all subsequent degree balances in
Lemma~\ref{vg:classical-regimes} therefore hold for $q_n(x)$.

In the Hermite case, $\Dd\mathbf u=R(x)\mathbf u$ makes $\op{H}$
symmetric. For all large $n$, degree preservation and orthogonality
make $\op{H}q_n(x)+(n/c)q_n(x)$ orthogonal to $\Pp_{n-1}$.
Its degree is at most $n-1$, so the nondegeneracy in
Lemma~\ref{lem:large-degree-duality} makes it zero. Thus
$\op{K}=\Lambda(-c\op{H})$: their difference annihilates monic
polynomials of every sufficiently large degree, which no nonzero
polynomial differential operator can do.

In the decreasing regime, the trace limit is $S(z)=(z-b)^{-1}$,
so the same degree comparison gives $f(x)=\kappa(x-b)^{2m}$.
The proof of Lemma~\ref{vg:classical-regimes}(iv) then uses only
the structural equation, the local pole bound, and the absence of a
polynomial annihilator from Lemma~\ref{lem:large-degree-duality}.
It therefore gives the Bessel conclusion for $q_n(x)$.

In the bounded and quadratic regimes, the dual representation
and absence of a polynomial annihilator in
Lemma~\ref{lem:general-foundation} also justify the proof of
Lemma~\ref{vf:classical-sources}.
The Hermite and Bessel regimes, the Laguerre regime with noninteger
residue, and the Jacobi regime with two noninteger residues give
classical $q_n(x)$ in all sufficiently large degrees.
In all remaining cases, that proof constructs
$$
 \mathbf v=\chi(x)\mathbf u,
$$
where $\chi(x)$ is the monic normalization of $f(x)$.
We retain this source even when the $q_n(x)$ themselves are classical.
Its classical operator has distinct eigenvalues in every degree.
Its monic eigenpolynomials $Q_n(x)$ are orthogonal
for $\mathbf v$. A zero norm would give
$\chi(x)Q_n(x)\mathbf u=0$, contradicting the absence of a polynomial
annihilator. Hence $\mathbf v$ is a regular Laguerre or Jacobi
functional with the parameter restriction of
Lemma~\ref{vf:classical-sources}.

Let $Q(x)=(Q_n(x))_{n\ge0}$ be the monic orthogonal sequence of
$\mathbf v$. Put $d=\deg\chi(x)$ and $k_n=\mathbf v(Q_n(x)^2)$.
Orthogonality gives the first expansion in
$$
 \begin{alignedat}{2}
 q_n(x)&=\sum_{j=0}^d a_j(n)Q_{n-j}(x),\quad &\chi(x)Q_n(x)&=\sum_{j=0}^d b_j(n)q_{n+j}(x).
 \end{alignedat}
$$
For the second expansion, use $\mathcal W$ and $q_j(x)$, $j\ge N$.
Pairings remove the terms below index $n$; the remaining term in
$\mathcal W$ is orthogonal to $\mathcal W$ and is zero by nondegeneracy. The extreme
coefficients satisfy $a_d(n)k_{n-d}=h_n$ and $b_0(n)h_n=k_n$, so
are nonzero. Comparison of fixed polynomial coefficients gives rational
$a_j(n),b_j(n)$.
The kernel argument in Lemma~\ref{vf:finite-connection} uses the dual
representation of Lemma~\ref{lem:general-foundation} in place of
\eqref{eq:regular-duality}: a functional vanishing on all sufficiently
large $q_n(x)$ is a polynomial multiple of $\mathbf u$.
Its uniqueness uses regularity of $\mathbf v$ and absence of a
polynomial annihilator of $\mathbf u$. The factorization proof applies
unchanged, since consecutive $q_n(x)$ are coprime. Denote the resulting connection by
$$
 \begin{gathered}
 q(x)=\op{A}_0Q(x),\\[7pt]
 \op{B}_0\op{A}_0=\chi(\op{J}_Q),\quad \op{A}_0\op{B}_0=\chi(\op{J}_q).
 \end{gathered}
$$
Here $q(x)=(q_n(x))_{n\ge N}$, and the identities are at sufficiently
large indices. Lemma~\ref{lem:bispectral-prefixes}, with the
orthogonality and duality just proved, gives polynomial differential
eigenoperators for all intermediate sequences in this connection.

\emph{Recovery of the common factor.}
The coefficient of $\partial^{2m-1}$ in
$\op{K}=G(x)^{-1}\op{L}_0G(x)$ gives
$$
 k_{2m-1}(x)=[\partial^{2m-1}]\op{L}_0
       +2mf(x)\frac{G'(x)}{G(x)}.
$$
Thus every root of $G(x)$ is a zero of $f(x)$. In the Hermite
case $f(x)$ is constant, so $G(x)=1$. In the Bessel case, write
$G(x)=(x-b)^h$ and $\op{K}=p(\op{H})$, where
$$
 \op{H}=t^2\partial^2+((\alpha+2)t+\eta)\partial,\quad
 t=x-b,\quad \eta\ne0.
$$
If $h>0$ and $\ell=\deg p(z)$, the lowest term of
$\op{H}^{\ell}(t^{-h})$ has coefficient
$\eta^{\ell}(-h)^{\underline\ell}\ne0$ and exponent $-h-\ell$.
It produces a pole in $\op{L}_0(1)=G(x)p(\op{H})(1/G(x))$, so $G(x)=1$.

At a Laguerre or Jacobi endpoint $a$ with noninteger residue
$\alpha$, $f(x)$ has a zero of order $m$. Write $t=x-a$ and let
$t^{-e}p(t\partial)$ be the part of $\op{K}$ giving the largest
decrease in the power of $t$. Then $e\ge m$, $\deg p(z)\le2m$, and
$$
 p(0)=\cdots=p(e-1)=0,\quad p(z)=p(e-\alpha-1-z).
$$
If $G(x)=t^h g(t)$, with $h\ge1$ and $g(0)\ne0$, polynomiality
of $G(x)\op{K}G(x)^{-1}$ also makes $p(z)$ vanish at
$-h,\ldots,e-1-h$. These two integer intervals contain at least
$e+1$ distinct roots. Their reflections are noninteger and disjoint
from them, giving $\deg p(z)\ge2e+2>2m$. Hence $G(x)$ has no root at
a noninteger endpoint. This also covers Jacobi with one integer and one noninteger residue.
Thus $G(x)=1$ in the Hermite and Bessel cases and whenever all
endpoint residues are noninteger, proving~(i) in these cases.
Otherwise every root of $G(x)$ is an endpoint with integer residue
for the source constructed above.

At each root $a$ of $G(x)$, write
$R(x)=\alpha_a/(x-a)+O(1)$.
Lemma~\ref{lem:common-factor-bound} gives
$$
 1\le g_a:=\ord_aG(x)\le\alpha_a.
$$
Such an $a$ is also a root of $\chi(x)$. Order the factorization
from $Q(x)$ to $q(x)$ so that $a$ is last, and let $r(x)$ be its
penultimate sequence. Lemma~\ref{lem:bispectral-prefixes} gives a
polynomial differential eigenoperator $\op{K}_a$ for $r(x)$.
The last upper factor $\op{V}_a=\op{T}+v_a(n)$ satisfies
$$
 (x-a)r_n(x)=q_{n+1}(x)+v_a(n)q_n(x).
$$
The rational function $v_a(n)$ is nonzero. If $q_n(a)=0$ at a large
index, this identity gives $q_{n+1}(a)=0$, contradicting
coprimality. Thus
$$
 \begin{alignedat}{2}
 A_n&:=\frac{q_{n+1}(a)}{q_n(a)}=-v_a(n),\quad &r_n(x)&=\frac{q_{n+1}(x)-A_nq_n(x)}{x-a}.
 \end{alignedat}
$$
Thus $\varphi_a(n)=q_n(a)\ne0$ eventually, with
$\varphi_a(n+1)/\varphi_a(n)=A_n\in\mathbb C(n)\setminus\{0\}$.

For $\mathbf u_a=(x-a)\mathbf u$, direct pairing gives
$$
 \begin{alignedat}{2}
 \mathbf u_a(r_n(x)p(x))&=0\quad\deg p(x)<n,\quad &\mathbf u_a(r_n(x)^2)&=-A_nh_n\ne0.
 \end{alignedat}
$$
This functional has no polynomial annihilator.
The degree and symmetry arguments of Lemma~\ref{lem:general-foundation}, using
Lemma~\ref{lem:large-degree-duality}, apply to $\op{K}_a$.
Its rational weighted transpose has
$$
 R_a(x)=R(x)+\frac1{x-a}.
$$
Indeed, the product rule gives
$(x-a)f(x)\Dd\mathbf u_a=f(x)(1+(x-a)R(x))\mathbf u_a$.
Comparing with the structural identity of $\op{K}_a$, and clearing
polynomial factors, proves the formula by the absence of an
annihilator.

Put $\op{J}=\op{J}_q$ and $\op{S}=\op{T}^{-1}$. Since
$(\op{J}-a)\varphi_a=0$, the functions
$$
 f_a(n)=-\frac{\varphi_a(n)}{\varphi_a(n-1)},\quad
 u_a(n)=\frac{C_n}{f_a(n)}
$$
give $\op{F}=1+f_a(n)\op{S}$, $\op{U}=\op{T}+u_a(n)$, and
$\op{J}-a=\op{U}\op{F}$. The two reversals are
$$
 \begin{gathered}
 \op{J}_1-a=\op{F}\op{U}
 =(\op{F}\op{T})(\op{S}\op{U}),\\[7pt]
 \op{F}\op{T}=\op{T}+f_a(n),\quad \op{S}\op{U}=1+u_a(n-1)\op{S}.
 \end{gathered}
$$
Both lower factors preserve monicity. The second reversal gives
$$
 \begin{alignedat}{2}
 \op{J}_2&=\op{S}\op{J}\op{T},\quad &(\op{S}\op{U})\op{F}q(x)&=(x-a)(q_{n-1}(x))_{n\ge N+1}.
 \end{alignedat}
$$
The new lower coefficients are $f_a(n)u_a(n-1)$ and
$C_{n-1}$, both nonzero rational functions. After $j$ pairs use the original solution $\varphi_a(n-j)$;
its nonvanishing avoids division by zero at repeated roots.
After $g$ pairs its $n$th polynomial is $G(x)q_{n-g}(x)=P_n(x)$,
of degree $n$.

Every intermediate also has a polynomial differential eigenoperator.
Let $H_j(x)$ be the product of the first $j$ restored factors,
with multiplicities. After $j$ complete pairs, its $n$th polynomial
is $H_j(x)q_{n-j}(x)$.
For $j>0$, Lemma~\ref{lem:polynomial-conjugates} applied to
$\op{K}$ gives its polynomial differential eigenoperator; for
$j=0$, use $\op{K}$ itself.

For the next endpoint $a$, use the sequence $r_n(x)$ just
constructed. The first reversal in the pair is precisely
$$
 H_j(x)\bigl(q_{n-j}(x)-A_{n-j-1}q_{n-j-1}(x)\bigr)
 =H_{j+1}(x)r_{n-j-1}(x).
$$
At $a$, the residue of $R_a(x)$ is $\alpha_a+1$; at every other
root of $G(x)$, it is unchanged. The multiplicities of $H_{j+1}(x)$
are within these bounds. Lemma~\ref{lem:polynomial-conjugates}
applied to $\op{K}_a,H_{j+1}(x)$ therefore gives a polynomial
differential eigenoperator for this sequence as well. The second
reversal returns to $H_{j+1}(x)q_{n-j-1}(x)$, already covered.
All eigenvalue polynomials remain nonconstant under the fixed index
shifts. Thus every step is bispectral, including repeated endpoint factors.

For the final sequence, put $d=\deg\chi(x)$ and define
$$
 \begin{gathered}
 \op{J}_P=\op{T}^{-g}\op{J}_q\op{T}^{g},\quad \op{A}_G=\op{T}^{-g}\op{A}_0G(\op{J}_Q),\\[7pt]
 \op{B}_G=\op{B}_0\op{T}^{g},\quad \Xi(x)=\chi(x)G(x).
 \end{gathered}
$$
The connection identities above give
$$
 \begin{gathered}
 P(x)=\op{A}_GQ(x),\quad \op{A}_G\op{J}_Q=\op{J}_P\op{A}_G,\\[7pt]
 \op{B}_G\op{A}_G=\Xi(\op{J}_Q),\quad \op{A}_G\op{B}_G=\Xi(\op{J}_P).
 \end{gathered}
$$
The extreme shifts of $\op{A}_G$ are $[-d-2g,0]$, with nonzero
coefficients. For Laguerre, apply the proof of
Lemma~\ref{vf:bispectral-composition} to these factors and the given
operator $\op{L}$. For Jacobi, $\chi(x)=\sigma(x)^m$, so $d=2m$;
use the centering shift $\op{T}^{m+g}$ in that proof.
It uses the displayed factor identities and $\op{L}$, without
coprimality of $P_n(x)$ or equal endpoint multiplicities of $\Xi(x)$.
It therefore gives the bispectral factorization with differential
operator $p(\op{L})$, where $p(z)$ is nonconstant, and recovers
$P(x)$ with its original degree labels.

As in the regular proof, degree comparison and subtraction of the
constant term put every intermediate eigenoperator in the
form~\eqref{eq:differential-operator-form}.
The number of steps and all parameters are fixed. Choose $N$ as in
the proof of Theorem~\ref{thm:regular-hst}, including all restoration
pairs. Then every reversal has rational triangular factors and a fixed
polynomial differential eigenoperator, and each resulting polynomial
is the unique monic eigenpolynomial in its original degree. Thus these
steps realize the same bispectral Darboux construction and prove~(i).
\end{proof}

\begin{proposition}\label{prop:general-canonical-exception}
With Bessel restricted to nonzero scale, the shifted canonical basis
$((x-b)^n)_{n\ge0}$ is neither a classical sequence nor obtainable
from a Laguerre or Jacobi sequence by finitely many rational Darboux
factor reversals with finite difference factors.
\end{proposition}

\begin{proof}
Equation~\eqref{eq:canonical-basis}, after translation, gives
$\op{J}=\op{T}+b\op{I}$ and minimum differential order one.
Every standard classical family at fixed parameters has a nonzero
rational lower recurrence coefficient for sufficiently large $n$.
An invertible affine change preserves this property, so none is the
shifted canonical basis.

For a nonzero finite rational difference operator $\op{A}$, let
$\ell(\op{A})$ be its least shift. The extreme coefficient in a product
cannot cancel, so
$$
 \ell(\op{A}\op{B})=\ell(\op{A})+\ell(\op{B})
 =\ell(\op{B}\op{A}).
$$
Thus a reversal $\op{J}-a=\op{U}\op{F}$,
$\widehat{\op{J}}-a=\op{F}\op{U}$ preserves least shift $-1$
whenever the lower recurrence coefficient is nonzero. Adding
$a\op{I}$ does not alter a negative least shift. Every operator in
a finite chain from a Laguerre or Jacobi source therefore has a nonzero
lower coefficient. It cannot equal $\op{T}+b\op{I}$.
\end{proof}

\section{The positive-definite case}
\label{sec:positive-classification}

A functional $\mathbf u$ is positive definite if it is real on
$\mathbb R[x]$ and $\mathbf u(p(x)^2)>0$ for every nonzero
$p(x)\in\mathbb R[x]$.

\begin{theorem}
\label{thm:positive-classification}
Let $\mathbf u$ be positive definite, with monic orthogonal sequence
$(P_n(x))_{n\ge0}$. There is a fixed nonscalar differential operator
$\op{L}$ with polynomial coefficients such that
$$
 \op{L}P_n(x)=\lambda_nP_n(x),\quad n\ge n_0,
$$
for some $n_0\ge0$, if and only if, after a real invertible affine
change of variable, $\mathbf u$ has one of the following forms,
for every $p(x)\in\Pp$:
\begin{enumerate}[label=(\roman*)]
\item Hermite:
\begin{equation}\label{eq:positive-hermite}
 \duality{\mathbf u}{p(x)}=c\int_{\mathbb R}p(x)e^{-x^2}dx,\quad c>0.
\end{equation}
\item Laguerre with a possible mass at zero:
\begin{equation}\label{eq:positive-laguerre}
 \duality{\mathbf u}{p(x)}
 =c\int_0^\infty p(x)x^\alpha e^{-x}dx+Mp(0),
\end{equation}
where $c>0$, $\alpha>-1$ and $M\ge0$. If $M>0$, then
$\alpha\in\mathbb Z_{\ge0}$.
\item Jacobi with possible endpoint masses:
\begin{equation}\label{eq:positive-jacobi}
  \duality{\mathbf u}{p(x)} =c\int_{-1}^1p(x)(1-x)^\alpha(1+x)^\beta dx +M_+p(1)+M_-p(-1),
\end{equation}
where $c>0$, $\alpha,\beta>-1$ and $M_+,M_-\ge0$.
If $M_+>0$, then $\alpha\in\mathbb Z_{\ge0}$; if $M_->0$, then
$\beta\in\mathbb Z_{\ge0}$.
\end{enumerate}
The representing positive measure is unique. Every such nonscalar
eigenoperator has even order $2m$. If $q$ is the density exponent at
an endpoint carrying a positive mass, then $2m\ge2q+4$.
\end{theorem}

\begin{lemma}\label{lem:positive-endpoint-mass}
Let $\mu$ consist of a positive Laguerre or Jacobi density $w(x)$,
with real exponents greater than $-1$, and nonnegative endpoint masses.
Let $\op{K}$ be a polynomial differential operator of order $2m>0$,
symmetric for $\mathbf u(p(x))=\int p(x)d\mu(x)$, with
$\op{K}(1)=0$ and $\op{K}^*=\op{K}$ for
$\partial^*=-\partial-w'(x)/w(x)$ on the open interval.
Suppose its leading coefficient has vanishing order $m$ at an
endpoint $a$ carrying a positive mass. If the density exponent there
is $q$, then $q\in\mathbb Z_{\ge0}$ and $2m\ge2q+4$.
\end{lemma}

\begin{proof}
Translate and reflect so that the endpoint is zero and the interval
lies to its right. Write $w(t)=t^qh(t)$, with $h$ analytic and
$h(0)>0$, and write $\op{K}=\sum_{j=1}^{2m}k_j(t)\partial_t^j$.
As in Lemma~\ref{vf:local-bounds}, let
$t^{-d}p(t\partial_t)$ be the part of $\op{K}$ giving the largest
degree decrease at zero. Then
$$
 \begin{gathered}
 d\ge m,\quad 0\ne p(z)\in\Pp_{2m},\\[7pt]
 p(0)=\cdots=p(d-1)=0,\quad p(z)=p(d-q-1-z).
 \end{gathered}
$$
Choose polynomials $U(t),V(t)$ such that
$$
 U(t)=t^d+O(t^{2m+1}),\quad V(t)=1+O(t^{2m+1}),
$$
and both vanishing through order $2m$ at the other finite endpoint.
Repeated integration by parts has boundary terms
$U^{(j-1-\ell)}(k_jwV)^{(\ell)}$, $0\le\ell<j$.
At zero their order is at least
$$
 q+d+\ord_0 k_j(t)-j+1\ge q+1>0;
$$
if the indicated derivative of $t^d$ is zero, its remainder is
$O(t^{q+2})$. At the other finite endpoint, with density exponent $\widetilde q>-1$
and distance $s$, the terms are $O(s^{\widetilde q+2m+3})$;
at an infinite Laguerre endpoint they vanish by exponential decay.
Polynomial interpolation gives such $U,V$. Thus the density
contribution to symmetry vanishes, as do the other atomic terms.
For $M=\mu(\{0\})>0$, symmetry gives
$$
 0=M\bigl((\op{K} U)(0)V(0)-U(0)(\op{K} V)(0)\bigr)
  =Mp(d).
$$
Reflection now gives the root sets
$\{0,\ldots,d\}$ and $\{d-q-1-j:0\le j\le d\}$.
If $q$ is noninteger or $q\ge d$, they are disjoint and give
$2d+2>2m$ roots. Hence $q\in\{0,\ldots,d-1\}$, and their union
is $[-q-1,d]\cap\mathbb Z$. Consequently
$$
 2m\ge\deg p(z)\ge d+q+2\ge m+q+2,
$$
which gives the asserted order bound.
\end{proof}

\begin{proof}[Proof of Theorem~\ref{thm:positive-classification}]
Positivity makes $\mathbf u$ regular and gives a positive Borel
representing measure $\mu$ of infinite support
\cite[Theorem~1.2]{shohattamarkin1943}. The polynomials $P_n(x)$
are real, $h_n>0$, $B_n\in\mathbb R$ and $C_n>0$.
Taking a real or imaginary part of the same positive order,
we may assume that $\op{L}$ is real.
Lemma~\ref{vf:foundation} gives a nonconstant eigenvalue polynomial. Use $\op{K},f,R$ from that
lemma and put $g(x)=k_{2m-1}(x)/m$.

Lemma~\ref{vg:recurrence-growth} gives
$|B_j|+\sqrt{C_{j+1}}\le A_0(j+1)$ for a fixed $A_0$.
Put $\|p(x)\|^2=\mathbf u(p(x)^2)$ for real polynomials.
In the orthonormal polynomial basis, multiplication by $x$ has
diagonal $B_j$ and off-diagonals $\sqrt{C_j}$. On real polynomials in $\Pp_k$, its
diagonal part and its two shift parts each have norm at most
$A_0(k+1)$. Thus, with $A=3A_0$,
$$
 \|xp(x)\|\le A(k+1)\|p(x)\|,\quad \deg p(x)\le k.
$$
Iteration, Cauchy--Schwarz and the exponential series give
\begin{equation}\label{eq:positive-moment-bound}
 \begin{gathered}
 \int x^{2n}d\mu(x)\le h_0 A^{2n}(n!)^2,\quad n\ge0,\quad \int |x|^n d\mu(x)\le h_0 A^n n!,\quad n\ge0,\\[7pt]
 \int e^{\varepsilon|x|}d\mu(x)\le\frac{h_0}{1-A\varepsilon},
 \quad 0<\varepsilon<1/A.
 \end{gathered}
\end{equation}
These bounds hold for every representing measure. Their Laplace
transforms are analytic on $|\Re z|<1/A$ and have the same Taylor
series at zero. The identity theorem makes them equal throughout
this strip; Fourier uniqueness then gives uniqueness of $\mu$.

The function
$$
 H(z)=-z\int f(x)e^{zx}d\mu(x)-\int g(x)e^{zx}d\mu(x)
$$
is analytic in a smaller strip. The functional identity gives
$H(0)=-\mathbf u(g(x))=0$ and, for $n\ge1$,
$$
 H^{(n)}(0)=-n\mathbf u(f(x)x^{n-1})-\mathbf u(g(x)x^n)=0.
$$
Thus $H=0$. Fourier inversion against smooth compactly supported
test functions gives the distributional identity below. The interchange
of integrals is justified by the finiteness of $f(x)\mu$ and $g(x)\mu$
and the rapid decrease of the Fourier transforms.
$$
 \Dd(f(x)\mu)=g(x)\mu.
$$
Away from the real zeros of $f(x)$, this reads $\Dd\mu=R(x)\mu$.
On each such interval choose $w>0$ with $w'/w=R$. Then
$\Dd(w^{-1}\mu)=0$, so $\mu=c_Iw(x)dx$ there, with $c_I\ge0$.
Only nonnegative atoms can remain at the finitely many real zeros
of $f(x)$.

In the bounded case of Lemma~\ref{vg:classical-regimes},
$b\in\mathbb R$ and $c>0$, so the two roots
$a=b-2\sqrt c$ and $d=b+2\sqrt c$ are real. The lemma gives
$$
 f(x)=\kappa((x-a)(x-d))^m,\quad
 R(x)=\frac{\rho_a}{x-a}+\frac{\rho_d}{x-d},
$$
with real residues. Bounded recurrence coefficients give a constant
$A_1>0$ such that $\int x^{2n}d\mu(x)\le h_0A_1^{2n}$ and
$\mu(\{|x|\ge r\})\le h_0(A_1/r)^{2n}$ for $r>A_1$.
Thus $\mu$ has compact support.
The densities on the exterior intervals must vanish. On $(a,d)$
the density is $c_0(x-a)^{\rho_a}(d-x)^{\rho_d}$, with $c_0>0$
because the support is infinite. Integrability gives
$\rho_a,\rho_d>-1$. Affine normalization yields
\eqref{eq:positive-jacobi}, apart from its integer conditions.

In the quadratic case, Lemmas~\ref{vg:classical-regimes} and
\ref{vf:classical-sources} give
$$
 f(x)=\kappa(x-a)^m,\quad
 R(x)=-\frac1\gamma+\frac\alpha{x-a}.
$$
Here $\gamma$ is real and nonzero, $\alpha$ is real, and the unique root $a$ of the
real polynomial $f(x)$ is real. On either side of $a$ the density
is a constant multiple of
$|x-a|^\alpha e^{-(x-a)/\gamma}$.
Finiteness eliminates the exponentially growing side; infinite
support makes the other constant positive. Integrability gives
$\alpha>-1$. The change $y=(x-a)/\gamma$ gives
\eqref{eq:positive-laguerre}, apart from its integer condition.

In the linear-growth case, $f(x)$ is constant and
$R(x)=-(x-b)/c$, with $b\in\mathbb R$ and $c>0$.
The measure equation gives a Gaussian density on the whole line,
hence \eqref{eq:positive-hermite} after normalization.
The decreasing case is impossible. Here
$f(x)=\kappa(x-b)^{2m}$ and
$R(x)=\alpha/(x-b)+\eta/(x-b)^2$, with real $b,\alpha,\eta$.
On each of $(-\infty,b)$ and $(b,\infty)$ the density is a
nonnegative constant multiple of
$|x-b|^\alpha e^{-\eta/(x-b)}$.
For $|x-b|\ge1$ the exponential is bounded below by
$e^{-|\eta|}$. Finiteness of all moments forces both constants to
vanish. Then $\mu$ is carried by $\{b\}$, contradicting infinite
support.

For any nonscalar eigenoperator of order $2m$, use its structural
pair $(f,g)$ from Lemma~\ref{vf:foundation}. Normalize an endpoint
carrying a mass $M>0$ to zero. Let $\mathbf w$ be the functional
defined by the absolutely continuous part, with density $w(x)$.
Put $\sigma(x)=x$ in the Laguerre case and
$\sigma(x)=x(L-x)$ in the Jacobi case on $(0,L)$.
In both cases $\Dd(\sigma(x)\mathbf w)=\tau(x)\mathbf w$
for a polynomial $\tau(x)$, and $\sigma(x)$ vanishes at every
finite endpoint. Thus
$$
 \Dd(\sigma(x)^2\mathbf u)
   =\sigma(x)(\tau(x)+\sigma'(x))\mathbf u.
$$
Comparing this with $\Dd(f(x)\mathbf u)=g(x)\mathbf u$,
regularity and cancellation of the nonzero polynomial $\sigma(x)$
give
$$
 \sigma(x)(g(x)-f'(x))=f(x)(\tau(x)-\sigma'(x)).
$$
Consequently $(f(x)w(x))'=g(x)w(x)$ on the open interval.
Choose $p(x)=x$ for Laguerre and $p(x)=x(1-x/L)^2$ for Jacobi.
Then $p(0)=0$, $p'(0)=1$, and both $p$ and $p'$ vanish at the
other finite endpoint. The density exponents exceed $-1$, so
$f(x)w(x)p(x)$ vanishes at both ends; at infinity use the Laguerre
exponential decay. Testing the structural equation on $p(x)$ gives
$$
 0=\mathbf u(f(x)p'(x)+g(x)p(x))
   =\int(f(x)w(x)p(x))'dx+Mf(0)=Mf(0).
$$
Hence $f(0)=0$. The shapes in Lemma~\ref{vg:classical-regimes}
give vanishing order exactly $m$.
Lemma~\ref{lem:positive-endpoint-mass} therefore gives the integer
conditions and order bounds for every nonscalar eigenoperator;
Lemma~\ref{vf:foundation} gives evenness.

Conversely, let $w(x)$ denote the classical weight in the corresponding
displayed formula. Then
$\mathbf u(p(x)^2)\ge c\int p(x)^2w(x)dx>0$ for
$0\ne p(x)\in\mathbb R[x]$, so all displayed functionals are positive definite.
With zero masses, the functional identity for the classical density
makes its second-order operator symmetric. Degree preservation and
orthogonality give the classical eigenrelations.

For a positive Laguerre mass write $\alpha=k\in\mathbb Z_{\ge0}$ and
$\mathbf u=c\mathbf L_k+M\delta_0$, where $\mathbf L_k$ is integration
against $x^ke^{-x}$ on $(0,\infty)$. Then
$x\mathbf u=c\mathbf L_{k+1}$ is a regular classical source.
Let $Q_n(x)$ be its monic orthogonal polynomials and put
$y_n=\mathbf u(Q_n(x))$. Their recurrence gives
$$
 y_{n+1}+(2n+k+2)y_n+n(n+k+1)y_{n-1}=0,\quad n\ge1.
$$
The two solutions $(-1)^nn!$ and $(-1)^nn!(n+1)_{k+1}$ are
independent. Their initial values give
$$
 y_n=(-1)^nn!\left(ck!+\frac{M}{(k+1)!}(n+1)_{k+1}\right),\quad n\ge0.
$$
In particular $y_n\ne0$. Orthogonality gives
$$
 P_n(x)=Q_n(x)+a_nQ_{n-1}(x),\quad
 a_n=-\frac{y_n}{y_{n-1}},\quad n\ge1.
$$
Thus $a_n$ and the first two lower coefficients of $P_n(x)$ are
rational in $n$. Lemma~\ref{vf:finite-connection} applies
with $\chi(x)=x$, and the Laguerre construction in the paragraph
following Lemma~\ref{lem:bispectral-prefixes}, which does not require a
previous eigenoperator of $P(x)$, gives the required polynomial
differential operator.

For Jacobi with at least one positive mass, interchange the endpoints
if necessary so that $M_->0$. Move to $(0,2)$ by $t=x+1$.
Retain $\mathbf u$ for the transformed functional and $P_n(t)$ for
the monic family obtained from $P_n(x)$ by substituting $x=t-1$.
Then
$$
 \mathbf u=c\mathbf J_{\beta,\alpha}+M_-\delta_0+M_+\delta_2,\quad
 t(t-2)\mathbf u=-c\mathbf J_{\beta+1,\alpha+1}=\mathbf v,
$$
where $\mathbf J_{\beta,\alpha}$ is integration against
$t^\beta(2-t)^\alpha$ on $(0,2)$.
Put $A=\beta+1$, $B=\alpha+1$ and $S=A+B+1$.
Then $A$ is a positive integer and $B>0$.
Let $Q_n(t)$ be the monic orthogonal polynomials of $\mathbf v$, with
recurrence coefficients $B_n,C_n$. Orthogonality gives
$$
 P_n(t)=Q_n(t)+a_nQ_{n-1}(t)+b_nQ_{n-2}(t),\quad n\ge2.
$$
Define
$$
 e_n=((t-2)\mathbf u)(Q_n(t)),\quad
 f_n=(t\mathbf u)(Q_n(t)).
$$
Applying these functionals to the recurrence gives
$$
 e_{n+1}+B_ne_n+C_ne_{n-1}=0,\quad
 f_{n+1}+(B_n-2)f_n+C_nf_{n-1}=0,\quad n\ge1.
$$
Let $G_e(0)=G_f(0)=1$ and define the nonzero hypergeometric terms by
$$
 \frac{G_e(n+1)}{G_e(n)}
 =-\frac{2(n+1)(n+B+1)}{(2n+S)(2n+S+1)},\quad
 \frac{G_f(n+1)}{G_f(n)}
 =\frac{2(n+1)(n+A+1)}{(2n+S)(2n+S+1)}.
$$
The first recurrence has the two independent solutions $G_e(n)$ and
$G_e(n)(n+1)_A(n+B+1)_A$, as substitution shows. Thus
$$
 e_n=G_e(n)p_e(n),\quad
 p_e(n)=c_0+c_1(n+1)_A(n+B+1)_A.
$$
If $M_+>0$, then $B$ is a positive integer. The analogous two independent
solutions of the second recurrence give
$$
 f_n=G_f(n)p_f(n),\quad
 p_f(n)=d_0+d_1(n+1)_B(n+A+1)_B.
$$
If $M_+=0$, the connection obtained by lowering the Jacobi parameter
$B$ by one gives $f_n=G_f(n)f_0$, so put $p_f(n)=f_0$.
Indeed, if its monic sequence is
$T_n(t)=Q_n(t)+\tau_nQ_{n-1}(t)$, comparison of subleading coefficients
gives $-\tau_n=G_f(n)/G_f(n-1)$, and $(t\mathbf u)(T_n(t))=0$ for
$n\ge1$. Both $p_e(n)$ and $p_f(n)$ are invariant under
$n\mapsto-n-S$. Neither is zero, since $e_0<0$ and $f_0>0$.

To use the Jacobi involution, first normalize at $t=0$:
$$
 \begin{gathered}
 Q_n(t)=\kappa(n)\psi_n(t),\quad
 \psi_n(t)=\sum_{j\ge0}
 \frac{(-n)_j(n+S)_j}{(A+1)_j j!}\left(\frac t2\right)^j,\\[7pt]
 \frac{\kappa(n+1)}{\kappa(n)}
 =-\frac{2(n+A+1)(n+S)}{(2n+S)(2n+S+1)}.
 \end{gathered}
$$
This is the normalization of Lemma~\ref{vf:bispectral-composition}
under $x=1-t$, with its source parameters $\alpha=A$ and $\beta=B$.
If
$\op{A}=1+a(n)\op{T}^{-1}+b(n)\op{T}^{-2}$ is the monic connection,
put $\widetilde{\op{A}}=\kappa^{-1}\op{A}\kappa$.
Solving
$e_n+a_ne_{n-1}+b_ne_{n-2}=0$
and $f_n+a_nf_{n-1}+b_nf_{n-2}=0$ gives
$$
 \widetilde{\op{A}}
 =1+\widetilde a(n)\op{T}^{-1}+\widetilde b(n)\op{T}^{-2},
$$
where
$$
 \begin{aligned}
 N(n)&=(n+B)(n+B-1)p_e(n)p_f(n-2)
       -(n+A)(n+A-1)p_e(n-2)p_f(n),\\[7pt]
 \Pi(n)&=(n+B-1)p_e(n-1)p_f(n-2)
       +(n+A-1)p_e(n-2)p_f(n-1),\\[7pt]
 \widetilde a(n)&=-\frac{nN(n)}{(n+A)(n+S-1)\Pi(n)},\\[7pt]
 \widetilde b(n)&=-\frac{n(n-1)(n+B-1)\Pi(n+1)}
 {(n+A)(n+S-1)(n+S-2)\Pi(n)}.
 \end{aligned}
$$
The leading coefficient of $\Pi(n)$ is twice that of $p_e(n)p_f(n)$,
so $\Pi(n)\ne0$. These formulas and the ratio of successive
$\kappa(n)$ show that $a_n,b_n$ and the first two lower coefficients of $P_n(t)$
are rational in $n$. Lemma~\ref{vf:finite-connection} therefore applies with
$\chi(t)=t(t-2)$.
The polynomial identities
$$
 N(-n-S)=-N(n+2),\quad \Pi(-n-S)=-\Pi(n+3)
$$
give
$$
 \widetilde a(-n-S)
 =\frac{\widetilde a(n+2)}{\widetilde b(n+2)},\quad
 \widetilde b(-n-S)=\frac1{\widetilde b(n+2)}.
$$
Consequently, for the involution of Lemma~\ref{vf:bispectral-composition},
$$
 \iota(\widetilde{\op{A}})
 =\widetilde b(n+2)^{-1}\op{T}^2\widetilde{\op{A}}.
$$
The even-width construction in the paragraph following
Lemma~\ref{lem:bispectral-prefixes} now gives a fixed nonscalar
polynomial differential eigenoperator for $P(t)$.

For comparison, the explicit Laguerre construction of
\cite{koekoeks1991}, in the form of
\cite[Proposition~2.3]{markett2024}, with $\nu=M/(ck!)$,
gives order $2k+4$ and leading coefficient
$(-1)^{k+1}\nu x^{k+2}/(k+2)!\ne0$.
For Jacobi, put
$$
 k^{-1}=\int_{-1}^1(1-x)^\alpha(1+x)^\beta dx,\quad
 m_\pm=kM_\pm/c.
$$
The construction in
\cite[Section~2, equations~(2), (4)--(19)]{koekoeks2000} gives,
with zero-mass terms omitted,
$$
 \op{L}=\op{L}_{\alpha,\beta}
       +m_-\op{A}+m_+\op{B}+m_-m_+\op{C},
$$
where $\op{A},\op{B},\op{C}$ are the correction differential operators
in the cited construction, and $\op{L}_{\alpha,\beta}$ is the Jacobi
differential operator
$(1-x^2)\partial^2+[\beta-\alpha-(\alpha+\beta+2)x]\partial$.
For the retained terms, equations~(14)--(19) give finite sums with
coefficients independent of $n$, nonzero leading coefficients and orders $2\beta+4$, $2\alpha+4$
and $2\alpha+2\beta+6$ for a mass at $-1$ only, at $1$ only
and at both endpoints, respectively. These orders need not be minimal.
Scaling the functional and monic normalization preserve the
eigenrelations, proving sufficiency.
\end{proof}

\section{Concluding remarks}
\label{sec:concluding-remarks}

A regular functional $\mathbf u$ is semiclassical if
$\Dd(\phi(x)\mathbf u)=\psi(x)\mathbf u$ for polynomials
$\phi(x),\psi(x)$, with $\phi(x)\ne0$. It is classical when such a pair
has $\deg\phi(x)\le2$ and $\deg\psi(x)\le1$
\cite[Definition~6.6]{castillo2026}. The structural identity
\eqref{eq:structural-equation} gives
$$
 \Dd(f(x)\mathbf u)=\frac{k_{2m-1}(x)}m\mathbf u,
 \quad f(x)=k_{2m}(x),
$$
so every regular functional in Theorem~\ref{thm:regular-hst} is
semiclassical. The converse fails even for positive functionals.
For $p(x)\in\Pp$, define the positive functional
$$
 \mathbf u(p(x))=\int_{\mathbb R}p(x)e^{-x^4}dx.
$$
It satisfies $\Dd\mathbf u=-4x^3\mathbf u$, whereas Lemma~\ref{vf:local-bounds} would force
$\deg f(x)+2m\le0$ for a nonscalar differential eigenoperator.

The recurrence assumption cannot be dropped. Define monic polynomials
$(R_n(x))_{n\ge0}$ by
$$
 F(x,t)=\sum_{n=0}^{\infty}R_n(x)\frac{t^n}{n!}
 =\exp\left(\frac{xt}{\sqrt{1-t^2}}\right).
$$
The identities $tF_t=x(\partial+\partial^3)F$ and
$xF=(1-t^2)^{3/2}F_t$ give, for $n\ge0$,
$$
 \begin{gathered}
 (x\partial+x\partial^3)R_n(x)=nR_n(x),\\[7pt]
 xR_n(x)=\sum_{j=0}^{\lfloor n/2\rfloor}
 (-1)^j\binom{3/2}{j}\frac{n!}{(n-2j)!}R_{n-2j+1}(x).
 \end{gathered}
$$
For $n=2m$, the coefficient of $R_1(x)$ is
$(-1)^m\binom{3/2}{m}(2m)!\ne0$.
Uniqueness of expansion in the monic basis therefore excludes any
fixed-length recurrence for multiplication by $x$, even in sufficiently
large degrees.

The minimum differential order is three, even in sufficiently large
degrees. Indeed, for $n\ge2$ the generating function gives
$[x^{n-1}]R_n(x)=0$ and
$[x^{n-2}]R_n(x)=n(n-1)(n-2)/2$.
A polynomial differential eigenoperator of order at most two on such
a tail preserves degree, by the degree comparison in
Lemma~\ref{vf:foundation}, and therefore has the form
$$
 \op{M}=(a_2x^2+b_2x+c_2)\partial^2+(a_1x+b_1)\partial+a_0.
$$
Comparison at degree $n-1$ gives
$n(n-1)b_2+nb_1=0$ for all sufficiently large $n$, so $b_2=b_1=0$.
Comparison at degree $n-2$ then gives
$$
 c_2=(n-2)\bigl((2n-3)a_2+a_1\bigr).
$$
This polynomial identity forces $a_2=a_1=c_2=0$.
Thus $\op{M}$ is scalar, whereas the displayed third-order operator
is nonscalar.

\begingroup
\renewcommand{\addcontentsline}[3]{}
\section*{Acknowledgements}
The first author wishes to thank P. Maroni and J. Petronilho for
teaching him everything he knows about the functional approach to
orthogonal polynomials. The second author thanks his coauthors
E. Horozov and M. Tater for their long and intensive work on the
project that resulted in \cite{hst2025}.

\section*{Declaration on the Use of Artificial Intelligence}

All mathematical ideas and results presented in this manuscript are the
authors’ own, developed through traditional mathematical research. Artificial
intelligence tools were used to improve the exposition,
independently check details of the proofs, and search for potential
counterexamples. The authors take full responsibility for the mathematical
content and correctness of the manuscript.
\endgroup

\end{document}